\documentclass[]{interact}
\usepackage[utf8]{inputenc}
\usepackage{epstopdf}% To incorporate .eps illustrations using PDFLaTeX, etc.
\usepackage{subfigure}% Support for small, `sub' figures and tables
\usepackage[spanish,es-tabla]{babel}
\usepackage{apacite}
\usepackage{multirow}
\usepackage{natbib}
\usepackage{color}
\usepackage[hidelinks,colorlinks=true,citecolor=blue]{hyperref}
\usepackage{color}
\theoremstyle{plain}% Theorem-like structures

\theoremstyle{definition}

\theoremstyle{remark}

\begin{document}

%\articletype{ARTICLE TEMPLATE}
\articletype{preprint}

\title{Aprendizaje de los números complejos usando diferentes sistemas de cálculo simbólico y un sistema de evaluación en línea en formación inicial de profesores}

\author{
\name{Jorge Gaona\textsuperscript{a}\thanks{CONTACTO Jorge Gaona. Email: jorge.gaona@upla.cl}, Silvia López\textsuperscript{b} y Elizanbeth Montoya Delgadillo\textsuperscript{c} }
\affil{\textsuperscript{a}Universidad de Playa Ancha, Valparaíso, Chile; \textsuperscript{b}Universidad de Viña del Mar, Viña del Mar, Chile; \textsuperscript{c}Pontificia Universidad Católica de Valparaíso, Valparaíso, Chile}
}

\maketitle

\begin{abstract}
En este artículo se estudió el trabajo matemático personal de 15 profesores en formación inicial, en primer año de una universidad pública en Chile, a partir de dos tareas sobre números complejos. Estas tareas se plantearon en un CAA (Moodlle - Wiris) y se propuso resolverlas utilizando distintos CAS, como GeoGebra, Symbolab, Photomath y Wolfram Alpha entre otros. Se observaron dificultades y potencialidades en el trabajo matemático de los estudiantes. En las dificultades, se observó que los estudiantes tenían problemas para interpretar la información de los CAS y del feedback del CAA debido a que los conocimientos matemáticos previos no eran lo suficientemente solidos para hacerlo. En las potencialidades, se observó que distintas características de la tarea, junto con la articulación de dos o más artefactos permitió a los estudiantes darle significado a los objetos matemáticos involucrados. 
\end{abstract}


\section{Introducción}

Una metáfora para ver el trabajo con tecnología en educación, y particularmente en matemáticas, es el mito de Sísifo \citep{Camus1942}. En este mito, el héroe, castigado por los dioses, se ve obligado a subir una roca hasta la cima de una montaña, roca que luego cae y la tiene que volver a subir indefinidamente. Los usuarios (profesores y estudiantes) de tecnología viven una situación análoga, porque continuamente, frente a la aparición nuevos \textit{artefactos digitales} o su actualización, deben aprender  o reaprender su uso. En este contexto, la problemática del trabajo matemático con \textit{artefactos digitales} se renueva constantemente \citep{Aldon-Tragalova2019,Flores2022}.

Desde un punto de vista teórico, en %Flores et al. 
\cite{Flores2022} se definen –desde una perspectiva instrumental, dentro de la teoría de los espacios de trabajo matemático \citep{Kuzniak2016c}– los \textit{artefactos digitales} para enseñar o aprender matemáticas como un conjunto de proposiciones caracterizadas por ser ejecutables por una máquina electrónica que posee una inteligencia histórica y una validez epistemológica relativa. En este trabajo nos centramos en tareas mediadas por dos clases de artefactos específicos y la relación entre ellos: los sistemas de cálculo simbólico (CAS, por sus siglas en inglés) y los sistemas de evaluaciones asistidas por computador o en línea (CAA, por sus siglas en inglés) que, a su vez, utilizan un sistema CAS para la validación de respuestas.

Los CAS han evolucionado y han sido ampliamente estudiados durante estos últimos 20 años. En términos de acceso, son un tipo de tecnología que en algunos países estaba ampliamente presente desde finales de la década del 2000 \citep[p. 243]{Trouche2000}. Mientras que, en otros países, más de una década después, su acceso era aún limitado \citep{OCDE2015}. Actualmente, existen muchos CAS que son accesibles desde un navegador web\footnote{Algunos ejemplos: \url{https://es.symbolab.com/}, \url{https://www.wolframalpha.com/}, \url{https://www.geogebra.org/calculator}, \url{https://calcme.com/a} o \url{https://photomath.com/} entre muchos otros. Estos CAS serán analizados con más detalle más adelante. }, sin necesitar siquiera instalación y menos una calculadora específica que los soporte. En términos técnicos, también han evolucionado los sistemas de entrada, permitiendo un ingreso de información más flexible, donde incluso, es posible ingresar, por ejemplo, una ecuación desde una fotografía \citep{Webel2016}. En cambio, otras características no han sufrido la misma evolución. Tomemos el caso de ecuaciones que por restricciones matemáticas o de programación solo pueden mostrarse soluciones aproximadas, incompletas o exactas según el CAS utilizado, manifestándose la relatividad de la validez epistemológica mencionada en %Flores et al. 
\cite{Flores2022}. 

Otra evolución importante es la información de salida que entregan estos sistemas y cómo se articulan con otros \textit{artefactos digitales}, en particular con los CAA. Si los trabajos de comienzos de siglo, los CAS entregaban principalmente un objeto de salida, como las soluciones de una ecuación o un gráfico \citep{Artigue2002,Hoyles2004,Kieran2006,Lagrange2000,Ruthven2002,Trouche2000,Yerushalmy1999}, actualmente, los sistemas han evolucionado entregando no solo la respuesta, sino que las soluciones paso a paso de un cálculo \citep{Barzel2019}, o en entornos interactivos para trabajar el razonamiento y la prueba automatizada \citep{Botana2015,Richard2019} o para automatizar la evaluación, validando respuestas equivalentes, diferenciando expresiones según algunas características (por ejemplo, si una expresión está factorizada o no) y dando un feedback en función del enunciado o de la respuesta del estudiante \citep{Gaona2018a,Gaona2021b,Gaona2021a,Sangwin2007,Sangwin2015}. 

Cada una de estas mejoras renueva las problemáticas asociadas al aprendizaje y la enseñanza con CAS, tanto por el valor epistémico y pragmático de su uso \citep{Artigue2002}, como por cuáles son los objetos matemáticos, las tareas y la organización de la clase que se involucran con los nuevos artefactos \citep{Pierce2010}. Al respecto, %Lagrange 
\cite{Lagrange2000} plantea que las técnicas de cálculo son necesarias para la conceptualización, pero las técnicas de papel y lápiz tienden a quedar obsoletas debido a la facilidad de uso de los comandos CAS. Esta obsolescencia es un problema porque las técnicas tradicionales ya no pueden desempeñar su papel en la conceptualización y las técnicas de ``botón'' no pueden asumir directamente este papel \cite[p. 117]{Lagrange2005b}. En esta misma línea, %Jankvist et al. 
\cite{Jankvist2019} plantean que el uso intensivo de comandos del tipo “solve” oculta una serie de procesos y fenómenos matemáticos diversos. Además, el proceso de razonamiento a priori que caracteriza el modo de trabajo matemático entrará en competencia con estrategias de ensayo y error más superficiales. 

En este contexto, parece ser que una pregunta pertinente sería ¿cómo aprender matemáticas a pesar del uso de la tecnología disponible para los estudiantes? En particular, cuando acceden a CAS disponibles en un navegador web. Más específicamente, ¿cuáles son las características de una tarea sobre un concepto matemático que permiten un trabajo matemático rico a partir del uso y articulación de diferentes \textit{artefactos digitales}?. A continuación se presenta una revisión de la literatura relacionada con el objeto matemático de estudio; el sistema de los números complejos. 

\section{\textit{Artefactos digitales} y números complejos}
Algunos estudios han usado \textit{artefactos digitales} \citep{Dazevedo2016,Caglayan2016} para promover la \textit{visualización} y la comprensión geométrica de los números complejos. En  %Caglayan 
\cite{Caglayan2016}, los resultados muestran que el uso de un entorno dinámico mejora la competencia matemática, en particular para pensar en las múltiples representaciones. En   %Dittman et al.
\cite{Dittman2017} se trabajó con estudiantes en formación inicial, para desarrollar el conocimiento geométrico sobre los números complejos y las funciones de variable compleja, para crear un juego utilizando un software de geometría dinámica, los resultados muestran que el uso del \textit{artefacto digital}, ayudó a la comprensión de funciones de variable compleja y la geométrica de los números complejos. En ambos trabajos, aunque se trabajó con tecnología, fue en entornos controlados con una situación diseñada ad-hoc y no específicamente con CAS.

De forma más general, a través de un trabajo matemático con distintas representaciones las personas asignan significados y pueden comprender las estructuras matemáticas, es por esto el interés didáctico de considerarlas \citep{radford1998signs}. Pero esta comprensión no es fácil para los estudiantes, ya que esto requiere niveles altos de abstracción y la confrontación de sus ideas construidas previamente en el sistema de los números reales \citep{Randolph2019}.

Considerando aspectos cognitivos del aprendizaje de los números complejos, y utilizando la teoría de los Modos de Pensamiento, el estudio de %Randolph y Parraguez
\cite{Randolph2019}, con estudiantes chilenos de nivel escolar y universitario, evidenciaron que, al plantearles tareas en lenguaje algebraico y gráfico, la mayoría de ellos privilegia los modos de pensamiento analítico-aritmético del sistema de los números complejos, y entienden los números complejos como un conjunto y no como un sistema numérico, en particular los estudiantes de nivel escolar, resolvieron las actividades desde los modos de pensar del sistema de los números reales y no desde el sistema de los números complejos. En el nivel universitario, los futuros profesores de matemáticas y licenciados en matemáticas, no alcanzaron completamente los otros dos modos de pensar (sintético-geométrico y analítico-estructural), plasmando la necesidad de continuar un trabajo matemático que potencie la articulación de estos distintos modos de pensar y, que dependen entre otros elementos, de trabajar distintas representaciones de los números complejos.  

En los números complejos se pueden utilizar cuatro formas de representación algebraica: imaginaria ($\sqrt{-1}=i$), cartesiana ($a+i \cdot b$), polar ($r \cdot (cos(\theta) + i \cdot sen(\theta)))$ y exponencial ($r \cdot e^{\theta \cdot i}$), las cuales pueden dar lugar a diferentes representaciones geométricas en un plano rectangular o polar. Hay que destacar que la conversión de la forma polar a la rectangular se puede hacer utilizando trigonometría del rectángulo, calculando la norma del número complejo y el ángulo que forma con el eje real. A la inversa, conociendo el ángulo y la norma del número complejo, se pueden obtener las coordenadas rectangulares calculando los catetos del triángulo que se forma entre el número complejo, el origen y la proyección del número complejo sobre el eje real (tomando en cuenta los signos dependiendo del cuadrante). Pero, para hacer la transformación a su forma de Euler se necesitan argumentos del análisis complejo, específicamente mediante las expansiones de Taylor para las funciones seno, coseno y exponencial aplicadas a una variable compleja se puede demostrar esta igualdad \cite[p. 23]{Artigue1992}. Esto implica que ciertas argumentaciones no serán posibles de desarrollar cuando los elementos del análisis no se hayan trabajado previamente. Finalmente, en la misma línea, el estudio de \cite{Karakok2015} con profesores de matemáticas de secundaria, exploró la concepción de las distintas formas de representación de los números complejos, y la manera en que transitan entre ellas. Usaron como referente teórico las concepciones objeto-proceso definidas por Sfard, los resultados muestran que estos profesores, trabajaron con la forma de Euler de los números complejos, pero no hubo evidencias de que tuvieran una concepción estructural de esta forma. Además, sus resultados indican que los profesores necesitan oportunidades para desarrollar una concepción dual (objeto-proceso) de cada forma.

Frente a la complejidad epistemológica del objeto y la ausencia de trabajos utilizando CAS, se diseña una tarea en un CAA, la pregunta de investigación que guía este trabajo entonces es la siguiente: ¿Cuál es el trabajo matemático de los estudiantes que se produce a partir de la articulación de diferentes \textit{artefactos digitales} en tareas sobre números complejos propuestas en un CAA?

En este estudio nos centraremos en la noción de número complejo a partir de la coordinación y confrontación de distintos \textit{artefactos digitales}, más específicamente nos interesa considerar las distintas representaciones (algebraicas y geométricas) que se encuentran presentes en el trabajo matemático del estudio de los números complejos, y que son movilizados por distintos \textit{artefactos digitales}, tales como; GeoGebra, Wolfram Alpha y Symbolab, entro otros. 

\section{Marco Teórico}
En primer lugar, nos referimos a la teoría del Espacio de Trabajo Matemático (subsección \ref{subsec_teoria_etm}), los elementos que la componen y que juegan un papel importante en nuestro estudio. A continuación, se conceptuailiza el concepto de  tarea en un CAA y se establece la relación con el ETM (subsección \ref{subsec_tareas_en_CAA}). 

\subsection{La teoría del ETM}\label{subsec_teoria_etm}
\begin{figure}
    \centering
    \includegraphics[width=0.9\textwidth]{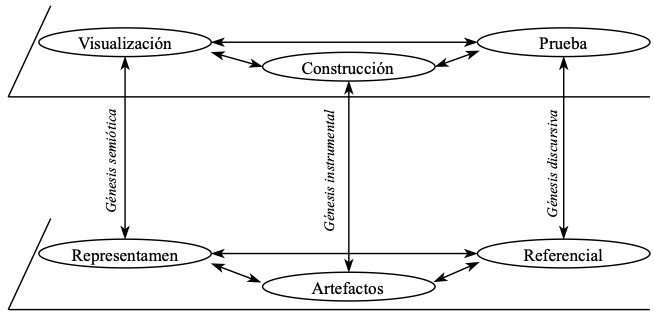}
    \caption{El modelo de los ETM extraído y traducido de \cite{Kuzniak2016c}}
    \label{fig:1_ETM}
\end{figure}

La teoría del ETM  \citep{KuzniakRichard2014,Kuzniak2016c,Kuzniak2022} nos permite describir, analizar, diseñar y comprender el trabajo matemático propuesto a un individuo en una institución y realizado por éste. Para caracterizar el trabajo matemático, la teoría considera dos planos: el  epistemológico y el cognitivo. En cada uno de estos planos se consideran tres componentes, organizados en tres \textit{génesis}, las cuales son esenciales en el trabajo matemático (ver figura \ref{fig:1_ETM}): 

\begin{itemize}
    \item La \textit{génesis semiótica}: que considera la \textit{visualización} como un proceso cognitivo donde un individuo da significado a los signos matemáticos (considerados, por ejemplo, en registros de representación semiótica) de la componente representamen.
    \item La \textit{génesis instrumental}: en la que un individuo transforma un artefacto (que puede ser material, digital o simbólico) en un instrumento para utilizarlo en un proceso cognitivo de construcción. En el caso de que sea digital, consideramos, a su vez, la inteligencia histórica y la validez epistemológica relativa con la que cuentan al presentar sus resultados \cite{Flores2022}.
    \item La \textit{génesis discursiva}: en la que un individuo utiliza el conocimiento matemático (definiciones, teoremas y propiedades) del \textit{referencial teórico} en el proceso cognitivo de las justificaciones.  
\end{itemize}

Estas tres \textit{génesis} están interrelacionadas. Se necesitan conocimientos para poder \textit{visualizar} los signos matemáticos o utilizar eficazmente un artefacto; la justificación utiliza representaciones semióticas de los objetos matemáticos; la construcción crea nuevos signos matemáticos, etc. La teoría del ETM pretende comprender el papel de cada una de estas \textit{génesis} y componentes en un sistema integrado que da lugar al trabajo matemático.  

La articulación y las interrelaciones entre las diferentes \textit{génesis} producen a través de \textit{planos verticales}. Cada uno de ellos destaca el papel preponderante de dos géneros: el plano semiótico-instrumental \textit{[Sem-Ins]}, el plano instrumental-discursivo \textit{[Ins-Dis]} y el plano semiótico-discursivo \textit{[Sem-Dis]}. 

Los artefactos desempeñan un papel crucial en nuestro estudio. Pueden ser de tres tipos y dependen del dominio matemático donde se trabaje:  
\begin{itemize}
    \item Materiales, como regla, compás o transportadores, que no se tratan aquí.
    \item Digitales, como las calculadoras simbólicas o gráficas como GeoGebra, Photomath, Symbolab y Wolfram Alpha o Calcme, entre otros. Estos últimos juegan un papel importante en nuestro estudio.
    \item Simbólicos, son rutinas o algoritmos que se utilizan en el proceso de construcción sin hacer ninguna conexión con la \textit{génesis discursiva} donde se \textit{justifican}.
\end{itemize}

En nuestra investigación los artefactos simbólicos los posicionamos, por ejemplo, como las reglas de factorización o la fórmula de segundo grado, entre otros. Éstos son muy útiles en el trabajo matemático, debido a su eficacia, sin embargo, es importante que el individuo sea capaz de reflexionar, de dar sentido a su uso. Por otro lado, los \textit{artefactos digitales} también juegan un papel preponderante en nuestro estudio, ya que permitirán obtener distintas soluciones a la ecuación propuesta y se podrá discutir sobre sus similitudes y diferencias.

\subsection{Tareas en un CAA}\label{subsec_tareas_en_CAA}

\begin{figure}[]
    \centering
    \includegraphics[height=0.9\textheight]{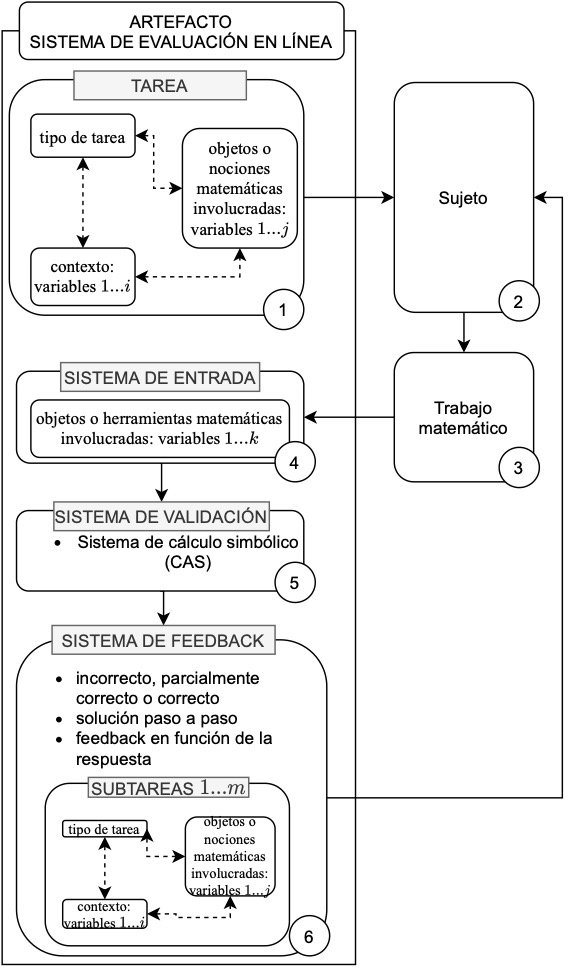}
    \caption{Articulación entre el artefacto CAA, los elementos didácticos de cada uno de sus componentes, el sujeto que se enfrenta a la tarea y el trabajo matemático que realiza}
    \label{fig:2_artefacto_tarea}
\end{figure}

En %Gaona 
\cite{Gaona2018b} se consideran tres componentes didácticas de una tarea: \textit{tipo de tarea}, \textit{objetos o nociones matemáticas} involucradas y \textit{contexto}. Por otra parte, en %Gaona 
\cite{Gaona2020} se hace una descomposición de un artefacto en un CAA en cuatro componentes: enunciado, sistema de entrada, sistema de validación y sistema de feedback, que están ligados al \textit{artefacto digital} involucrado. Al complementar estas dos descomposiciones más el sujeto que responde una tarea, quien a partir de esa tarea realiza un trabajo matemático específico (sección \ref{subsec_teoria_etm}), se pueden articular en el esquema que se muestra en la figura \ref{fig:2_artefacto_tarea}. 

En la tarea (ver 1 de la figura \ref{fig:2_artefacto_tarea}), se distinguen tres componentes, el \textit{tipo de tarea}, los \textit{objetos o nociones matemáticas} y el \textit{contexto}. 

El \textit{tipo de tarea} se identifica por un verbo y una acción precisa \citep{Chevallarda}, en nuestro caso, encontrar las raíces de una ecuación. El segundo componente, es el o los \textit{objetos o nociones matemáticos} involucradas que permiten precisar la tarea. En el caso de una ecuación, podrían ser: polinomial o trigonométrica y dentro de la ecuación hay una serie de variables que se pueden definir. Este objeto o herramienta es observable gracias a su representación semiótica \citep{Duval1995} y los elementos que la componen. Por ejemplo, si tenemos una ecuación polinomial de grado 2, la ecuación podría ser representada de forma gráfica, en lenguaje natural o algebraicamente o incluso en una representación dinámica a partir de un applet de geometría dinámica y los coeficientes podrían ser enteros, racionales, irracionales, entre otros. Cada una de esas elecciones puede influir en la naturaleza de las soluciones y en la dificultad de la tarea. Finalmente, el \textit{contexto} comprende los elementos o variables que permiten dar significado a los objetos o nociones involucradas. Este aparece conectado con una línea punteada, ya que no siempre se utiliza en una tarea, pero en caso de que exista, este puede influir en el tipo de tarea o en el objeto/herramienta matemática puesto en juego. A modo de ejemplo, si la variable de la ecuación representa una distancia entonces solo serán válidas las soluciones positivas.

A partir de esta tarea, hay un sujeto (ver 2 y 3 de la figura \ref{fig:2_artefacto_tarea}) que realiza un trabajo matemático. Este trabajo puede realizarse en interacción con el CAA, como también fuera de él. La descomposición conceptual de este trabajo es el que se realizó en el punto \ref{subsec_teoria_etm}.  

Luego, el sujeto ingresa una respuesta mediante el sistema entrada de la tarea (ver 4 en figura \ref{fig:2_artefacto_tarea}). En %Gaona 
\cite{Gaona2020} se identifican varios de estos de sistemas de entrada: espacio en blanco para ingresar una respuesta, editor de ecuaciones o sistema de reconocimiento de escritura a mano alzada para ingresar las respuestas, sistema “arrastrar y soltar”, sistema geométrico para responder o respuesta de opción múltiple. Del sistema de entrada dependerá el o los objetos matemáticos que se esperan como respuesta y su representación.  

Una vez ingresada la respuesta, el sistema realiza el proceso de validación (ver 5 en figura \ref{fig:2_artefacto_tarea}). Salvo el caso de selección múltiple, los CAA usados en matemáticas, deben contar con sistemas de cálculo simbólico o geométrico para comparar la respuesta ingresada con la que está definida en el sistema como correcta. 

Una vez que el sistema valida, puede entregar un feedback a los estudiantes (ver 6 en figura \ref{fig:2_artefacto_tarea}). Según %Hattie 
\cite{Hattie2008} los tipos de feedback en una interacción educativa son variados, en cambio, en un CAA son más limitados. Sin embargo, al menos se pueden identificar un feedback que le informa al estudiante si la respuesta es correcta, parcialmente correcta o incorrecta, dar algún comentario de tipo socioemocional, dar una solución paso a paso del problema planteado o dar un feedback en función de la respuesta del estudiante. El primer tipo de feedback, junto con los dos últimos son de particular interés, porque pueden ayudar o dificultar la creación de significados matemáticos. Este feedback lo recibe el sujeto y se confrontará al trabajo matemático que realizó, produciéndose similitudes y diferencias entre lo que hizo y lo que se esperaba que hiciera. También puede ser formativo, modificando el trabajo matemático inicial \citep{Gaona2021b} o podría producir dificultades puesto que el feedback entregado puede ser matemáticamente inadaptado al trabajo realizado por los sujetos, adaptado, pero mal interpretado o adaptado pero insuficiente \citep{Cazes2006}. El feedback más que generar certezas en los estudiantes, les ayuda a tener de aquello que no saben. En este feedback también se pueden identificar subtareas que, a su vez, se descomponen en \textit{tipo de tarea}, \textit{objetos o nociones matemáticas} y el \textit{contexto}. 

\section{Elementos metodológicos y contexto de la experimentación}

La investigación es de enfoque cualitativo lo cual será abordado mediante un estudio de casos, donde cada caso corresponde a un caso instrumental \citep{Bikner-Ahsbahs,yin2009case}. Se intervino en el ETM personal de los estudiantes en un curso de tecnología. Las unidades de análisis fueron 14 estudiantes de primer año en formación inicial de profesores, de una universidad pública chilena. La experimentación se desarrolló durante el primer semestre del 2021, en tres sesiones de clases de 135 min cada una, en modalidad online, sincrónica y vía plataforma Zoom\footnote{https://zoom.us/}.  Cabe señalar, que uno de los investigadores fue el profesor que implementó la situación de aprendizaje.  

Se plantea una situación didáctica compuesta por dos tareas que se presentaron en una plataforma Moodle/Wiris que se describen en la sección \ref{Subsec_Las_tareas}. A los estudiantes se les propuso trabajar con diferentes CAS (GeoGebra, Symbolab y Wolphram Alpha, entre otros).

 Se transcribieron todas las videograbaciones y los diálogos, tanto de los estudiantes como del profesor. Se utilizaron nombres de fantasía para los estudiantes, con el propósito de guardar su identidad. Las trascripciones de los diálogos referidas a expresiones matemáticas, se agregaron expresiones algebraica entre corchetes cuadrados, con el propósito de comunicar lo que sujeto dijo matemáticamente al oral, por ejemplo; “escribí $z$ elevado a tres igual a menos dos [$z^3=-2$]”
 
 %Por último, se añaden acciones hechas sobre el ordenador, tanto por estudiantes como por el profesor, mostradas como capturas de pantalla, con el propósito de complementar el análisis del trabajo matemático realizado por los estudiantes.
 
Las videograbaciones se subdividieron en episodios, para identificarlos consideramos dos criterios: 

\begin{itemize}
    \item Intervención de un estudiante que genere una discusión del trabajo matemático 
    \item Intervención del profesor (para un cambio de trabajo de la tarea y cambio del trabajo de los estudiantes) 
\end{itemize}

A partir de la teoría del ETM se levantaron categorías de análisis considerando aspectos cognitivos y epistemológicos del trabajo matemático que son considerados durante la implementación de la situación de aprendizaje.

\subsection{La tabla de análisis de los episodios}
Para el análisis de los episodios, se consideran los elementos del ETM que se encuentran presentes en el trabajo matemático. Las \textit{génesis} que son activadas: \textit{semiótica}, \textit{instrumental} y \textit{discursiva}, y los planos movilizados: \textit{[Sem-Ins]}, \textit{[Sem,Dis]} y \textit{[Ins-Dis]}. Los cuales fueron definidos en la tabla \ref{tab:1_ Categ_de_analisis}.

\begin{table}[h]
\caption{Categorías de análisis} 
\label{tab:1_ Categ_de_analisis}
\begin{center}
\scalebox{0.65}{\begin{tabular}{p{3cm}p{3cm}p{12cm}}
\toprule
\textbf{Categoría}     & \textbf{Componente}               & \textbf{Descriptor}                                                                                                                                                                                                                                                                                                                                                                            \\ \midrule
\textit{Génesis}  \newline \textit{semiótica}    & \textit{Representamen}                     & Considera   los elementos matemáticos en sus distintas formas de representación.                                                                                                                                                                                                                                                                                                               \\ 
   & \textit{visualización}                     & Los   objetos matemáticos son interpretados mediante procesos cognitivos ligados a   la \textit{visualización}, ya sean, gráfico, algebraico o lenguaje natural.                                                                                                                                                                                                                                        \\ \midrule
\textit{Génesis}  \newline \textit{instrumental} & \textit{Artefactos}                        & Se   consideran \textit{artefactos digitales} como programas informáticos, tales como;   GeoGebra, Symbolab, Wolfram Alpha, Calcme y Photomath.   Además, artefactos simbólicos, como algoritmos que   se utilizan en el proceso de construcción sin hacer ninguna conexión con la   \textit{génesis discursiva} donde se \textit{justifican}.                                                                          \\ 
 & \textit{Construcción}                      & Consideradas   como acciones que son hechas a través del uso de los distintos artefactos y   sus técnicas de uso asociadas.                                                                                                                                                                                                                                                                    \\ \midrule
\textit{Génesis} \textit{discursiva}   & \textit{Referencial} \newline \textit{teórico}                       & Asociados   a las definiciones, propiedades o teoremas   matemáticos pre-establecidas.                                                                                                                                                                                                                                                                                                           \\ 
   & \textit{Prueba}                            & Es   necesaria una argumentación, justificación o demostración para el razonamiento \textit{discursivo}, ya sea pragmático con el uso de algún artefacto   digital, o intelectual por medio de un razonamiento matemático (propiedad o   teorema).                                                                                                                                                    \\ \midrule
\textit{{[}Sem-Ins{]}}          & \textit{génesis} \textit{semiótica} e \textit{instrumental}  & Estas génesis se verán articuladas, en el caso por ejemplo; los objetos matemáticos son   representados por alguno de los \textit{artefactos digitales} usados en la tarea   (GeoGebra, , Symbolab, Wolfram Alpha, Calcme y Photomath), lo cual permite su   \textit{visualización} en distintos registros de representación semiótica (algebraico   y/o gráfico).                                                                                                           \\ \midrule
\textit{{[}Sem-Dis{]}}          & \textit{génesis} \textit{semiótica} y \textit{discursiva}    & Estas génesis se verán articuladas, en el caso por ejemplo; la representación de las soluciones de la ecuación obtenidas por los artefactos digitales o simbólicos en los registros algebraico y gráfico, son utilizados para   justificar algún resultado.                                                                                                                                                                                             \\ \midrule
\textit{{[}Ins-Dis{]}}          & \textit{génesis} \textit{instrumental} y \textit{discursiva} & Estas génesis se verán articuladas, en el caso por ejemplo; los \textit{artefactos digitales} dan información en los registros algebraico y gráfico   sobre las soluciones de la ecuación, esta información es usada, y en   coordinación con los conocimientos previos es posible dar una respuesta. Por   otro lado, es posible usar artefactos simbólicos como, por ejemplo, reglas de   factorización, para validar o justificar lo obtenido por los \textit{artefactos digitales}. \\ \bottomrule
\end{tabular}}
\end{center}
\end{table}

\subsection{Las tareas}\label{Subsec_Las_tareas}

\begin{figure}
    \centering
    \includegraphics[width=0.99\textwidth]{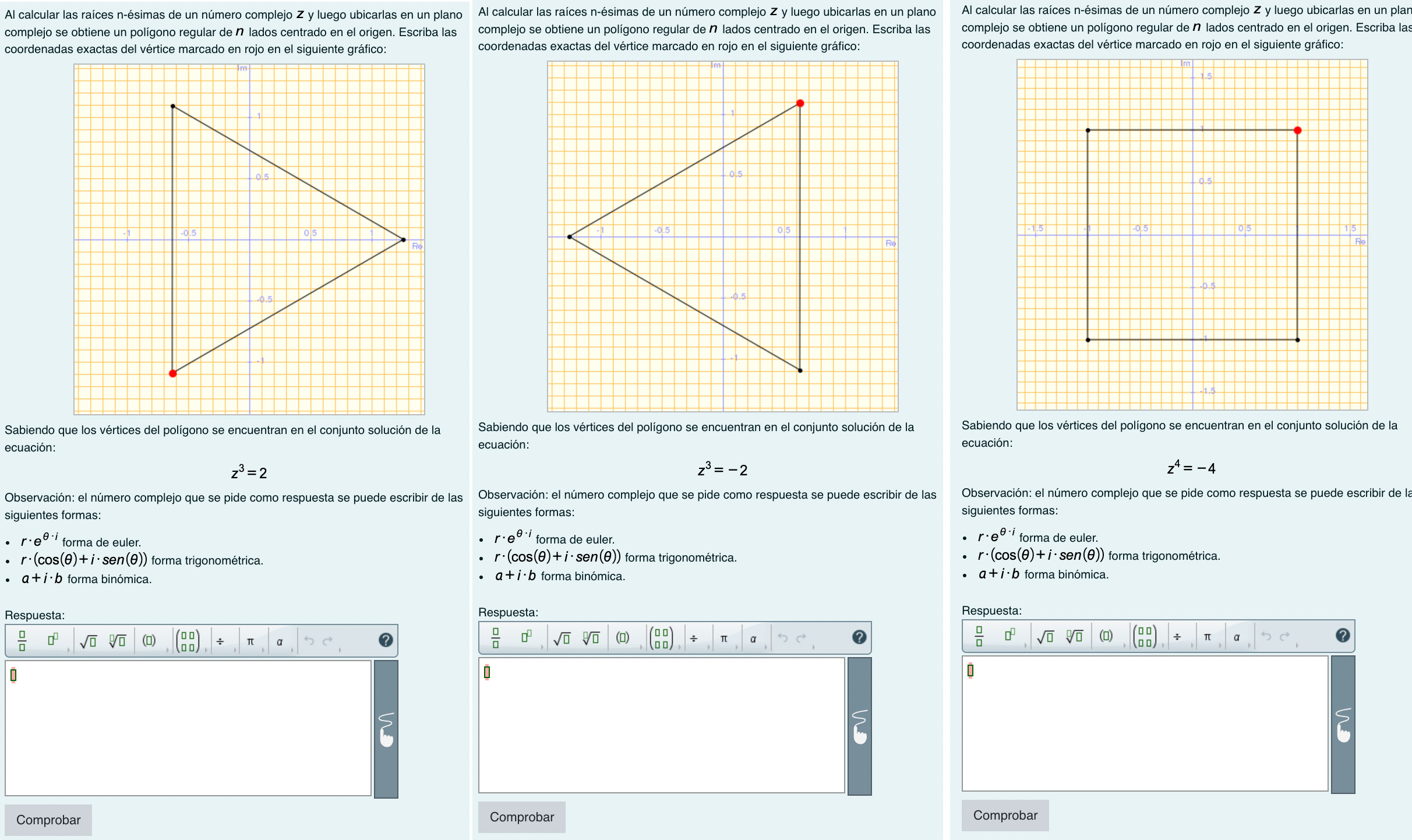}
    \caption{Enunciados de la tarea presentada a los estudiantes en Moodle/Wiris. En el de la izquierda, se pide la solución del tercer cuadrante de la ecuación $z^3=2$. En el del centro, se pide la solución del segundo cuadrante de la ecuación $z^3=-2$ y en el de la derecha, se pide la solución del primer cuadrante de la ecuación $z^4=-4$.}
    \label{fig:3_enunciados}
\end{figure}

El trabajo propuesto se desarrolló en torno a dos tareas principalmente, que constituyen la situación de aprendizaje. Cabe señalar, antes de la implementación hubo una primera parte de diagnóstico, dónde el profesor determinó el estado actual de conocimientos que tenían los estudiantes de los números complejos. A continuación se definen las tareas:
\begin{itemize}
    \item \textbf{Tarea 1:} Escribir, en el espacio dispuesto para ello, el número complejo marcado en rojo de la ecuación $z^3=2$ (figura \ref{fig:3_enunciados} izquierda) o $z^3=-2$ (figura \ref{fig:3_enunciados} del centro). Los valores $2$ o $-2$  al que está igualado la ecuación, al igual que el punto marcado en rojo son definidos mediante un algoritmo de forma aleatoria. La raíz marcada que aparece en el plano tiene como condición que la parte imaginaria debe ser no nula, por lo que nunca estará en el eje real.
    \item \textbf{Tarea 2:}  Escribir, en el espacio dispuesto para ello, el número complejo marcado en rojo de la ecuación $z^n=a$, donde $n$ puede tomar el valor 3 o 4 (figura \ref{fig:3_enunciados} derecha) y a podía tomar un valor entero entre $-9$ y $9$, incluyendo los extremos y excluyendo el cero. 
\end{itemize}

\begin{figure}
    \centering
    \includegraphics[width=0.8\textwidth]{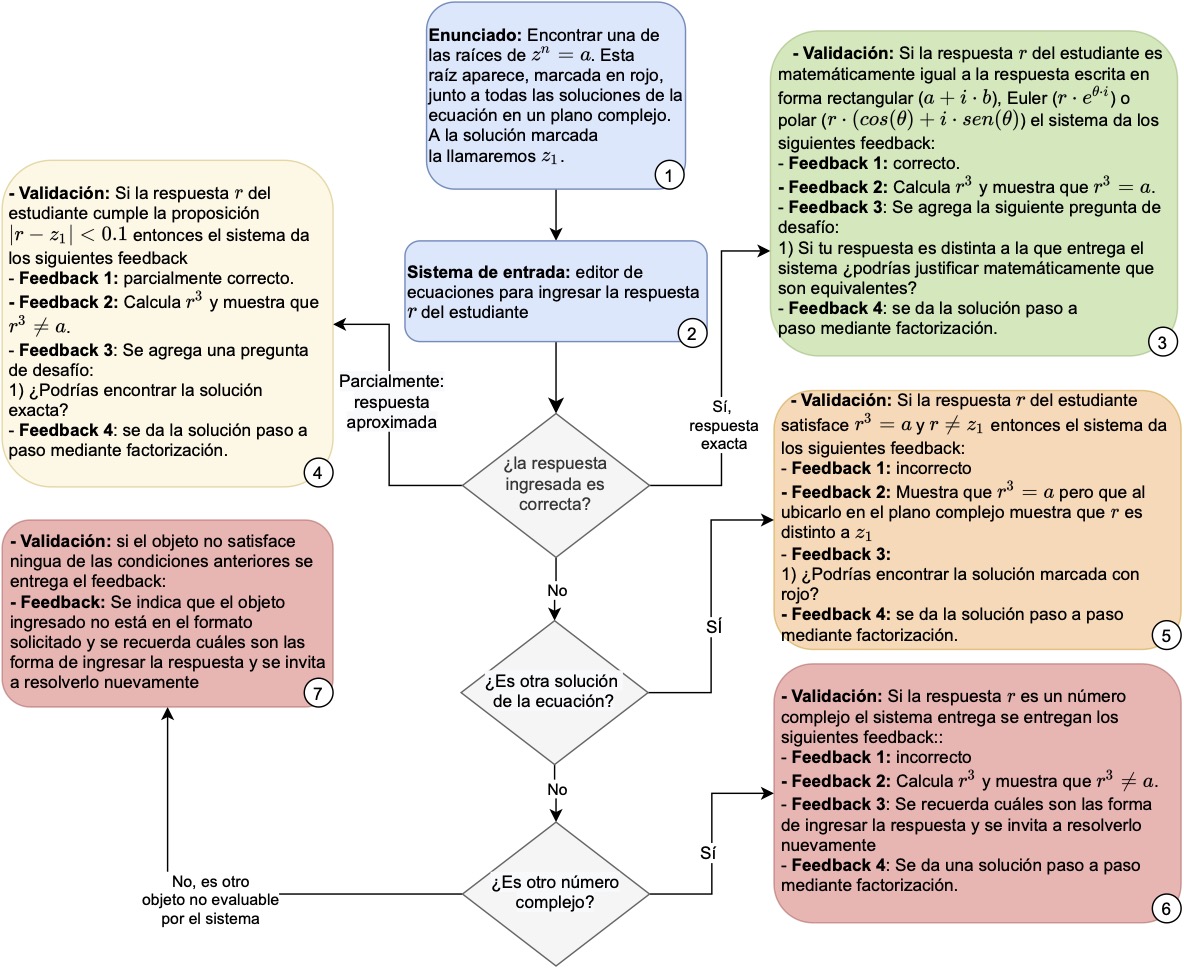}
    \caption{Árbol de decisiones para evaluar y dar un feedback a partir de las respuestas de los estudiantes. En 1 está en el enunciado, en 2 el sistema de entrada. En 3 la validación y feedback de la respuesta correcta, en 4 las respuestas parcialmente correctas y en 5, 6 y 7 de las respuestas incorrectas}
    \label{fig:4_arbol_de_decisiones}
\end{figure}

\newpage
En las tareas 1 y 2, una vez que el estudiante ingresa una respuesta al espacio dispuesto para ello, el sistema le entrega un feedback diferente dependiendo si la respuesta es correcta, parcialmente correcta o incorrecta. Estas diferentes opciones se esquematizan en el árbol de decisiones de la figura
\ref{fig:4_arbol_de_decisiones}.

\begin{figure}
    \centering
    \includegraphics[width=0.8\textwidth]{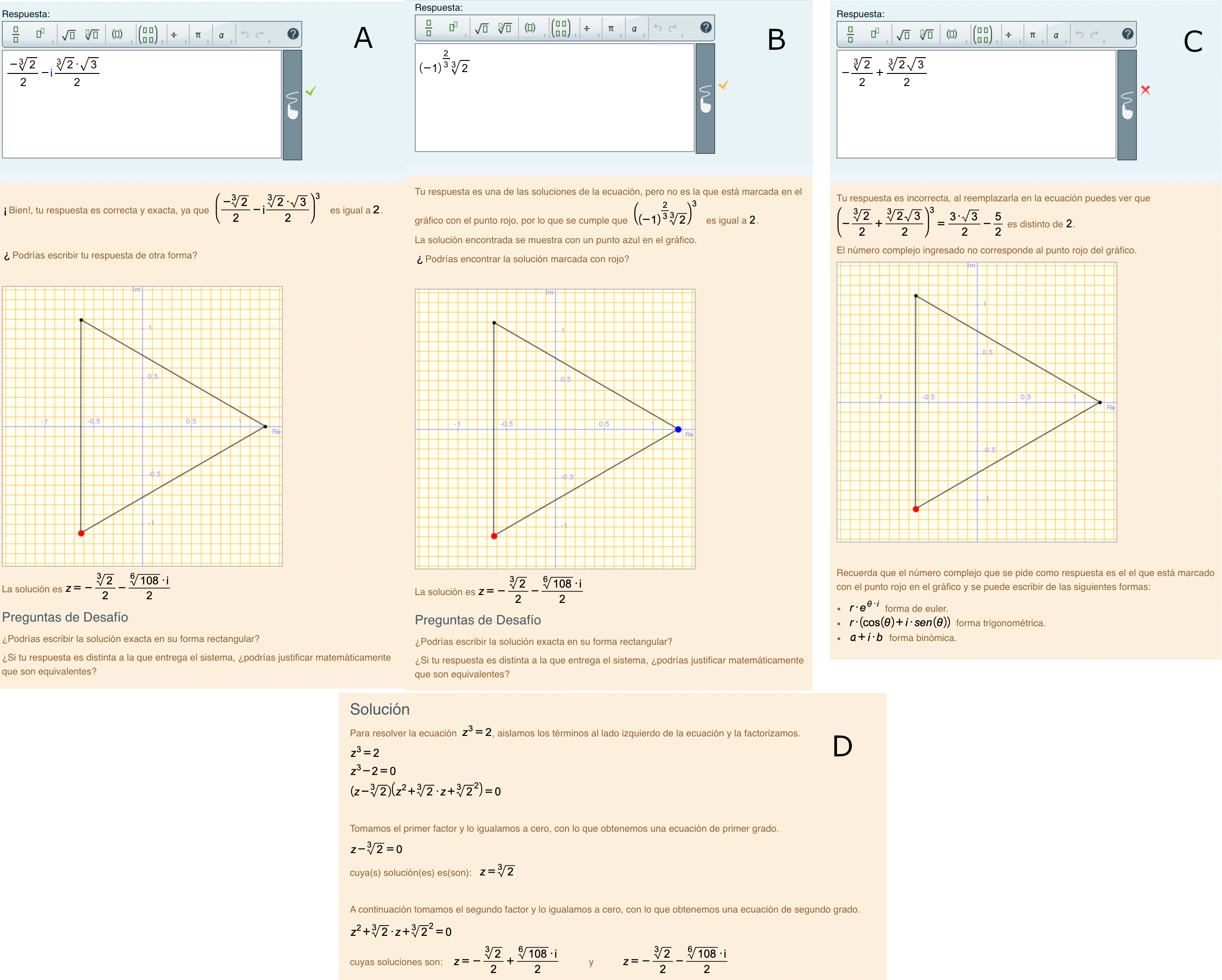}
    \caption{Validación y feedback para 3 casos. En A, la respuesta está escrita de forma rectangular y es correcta. En B, la respuesta ingresada es la solución real de la ecuación, pero no la que está marcada en rojo. En C, la respuesta es incorrecta y no es ninguno de los casos anteriores y en D, es un feedback que aparece en todas las versiones.}
    \label{fig:5_feedback_plataforma}
\end{figure}

A modo de ejemplo, en la figura \ref{fig:5_feedback_plataforma}, se muestran la validación y el feedback desplegados para los casos 4 (figura \ref{fig:5_feedback_plataforma}A) y 6 (figura \ref{fig:5_feedback_plataforma}B y C) del árbol de decisiones de la figura \ref{fig:4_arbol_de_decisiones}.  Al descomponer el artefacto tarea asociados a las tareas 1 y 2 en sus distintos componentes, se puede observar que:

\begin{itemize}
    \item \textbf{El enunciado:}  está compuesto por elementos aleatorios: en la ecuación $z^n=a$, donde $a$ toma los valores $-2$ o 2 en la tarea 1 y los números enteros entre $-9$ y 9 en la tarea 2. El exponente $n$ en cambio, en la pregunta 1 siempre vale 3 y en la tarea 2 varían entre $n=3$ o $n=4$. 
    \item \textbf{El sistema de entrada:} es un editor de ecuaciones para escribir los números complejos de forma rectangular: $a+i \cdot b$, polar: $r \cdot (cos(\theta)+i \cdot sen(\theta)) $ o Euler: $r \cdot e^{\theta*i}$.
    \item \textbf{El sistema de validación:} evalúa respuestas matemáticamente equivalentes, diferencia entre una respuesta aproximada y una exacta, y diferencia entre las soluciones de la ecuación ingresadas.
    \item \textbf{El sistema de feedback:} entrega información en función de la respuesta, y también un paso a paso en función de los parámetros de la tarea.
\end{itemize}

.

\subsection{Análisis a priori de los \textit{artefactos digitales}}\label{Subsec_Artefactos_digitales}

A los estudiantes se les propuso utilizar diferentes \textit{artefactos digitales} a disposición para resolver la tarea. Podían usar los que tuvieran a disposición y, además, se les entregó el enlace o se les propuso que instalaran en su celular alguna de estas aplicaciones para aproximarse a la solución del problema. 

\begin{itemize}
    \item \textbf{Photomath:} es una aplicación que permite resolver ejercicios como ecuaciones, operaciones y cálculos a partir de una foto del problema. Para este caso, al sacar una foto de $z^3=2$ da como resultado la solución real de forma exacta $\sqrt[3]{2}$ y de forma aproximada $1.259$.
    \item \textbf{Wolfram Alpha:}  es una aplicación web que funciona con Mathematica como motor de cálculo. Al escribir $z^3=2$, la página entrega información en distintos registros semióticos. Muestra\footnote{El software actualizó la información que muestra. Si ahora se ingresa $z^3=2$ muestra las soluciones complejas de forma correcta, además, agregaron un botón que permite elegir entre forma de Euler o Polar, sin embargo, al momento que se tomaron los datos, la información que aparece es la que se describe en este punto.} una forma alternativa de la misma ecuación: $z^3-2=0$, la solución real $\sqrt[3]{2}$  en la recta real (de forma gráfica), las soluciones complejas, que son erróneas porque son las mismas soluciones reales escritas de otra forma: $\sqrt[3]{-2}$ y $(-1)^{2/3}\sqrt[3]{2}$. Finalmente, muestra la ubicación de las soluciones en el plano complejo. Al presionar el botón ``Aproximate form'' se despliegan las aproximaciones de las soluciones en su forma rectangular. El resto de la información es la misma.
    
    También se puede hacer un click en alguna de las soluciones complejas y Wolfram Alpha 	entrega la forma de Euler, la forma polar y una forma rectangular alternativa. Las tres expresiones son correctas, sin embargo, en esta misma página sigue apareciendo la solución compleja como real. 
    \item \textbf{Symbolab:} Esta aplicación está disponible para móviles y en una versión web. Al ingresar $z^3=2$ muestra las tres soluciones en su forma algebraica, rectangular y exacta, las muestra factorizada por la raíz cúbica de 2:  
    $z=\sqrt[3]{2}$, $z=\sqrt[3]{2}\dfrac{-1+\sqrt{3}i}{2}$, $z=\sqrt[3]{2}\dfrac{-1-\sqrt{3}i}{2}$ 
    y con la multiplicación desarrollada: 
    $z=\sqrt[3]{2}$, $z=-\dfrac{\sqrt[3]{2} }{2}+i \dfrac{\sqrt[3]{2} \cdot \sqrt{2}}{2}$, $z=-\dfrac{\sqrt[3]{2} }{2}-i \dfrac{\sqrt[3]{2} \cdot \sqrt{2}}{2}$.
    \item \textbf{GeoGebra:}  es un \textit{artefacto digital} que tiene una versión de escritorio y otra versión web. En ambas, el comando \textit{Resuelve($z^3=2$)} entrega la solución en los reales: $\sqrt[3]{2}$. En cambio, el comando \textit{RaízCompleja($z^3-2$)} entrega todas las soluciones de la ecuación de forma aproximada: $z_1=-0.63-1.09i$, $z_1=-0.63+1.09i$, $z_1=-1.26+0i$ y también, muestra la ubicación de estas en el plano cartesiano. Hay otro comando que se llama \textit{SolucionesC($z^3=2$)} y da como soluciones: $\left\lbrace \dfrac{1}{2}(i \sqrt[3]{2}\sqrt{3}-\sqrt[3]{2}),\dfrac{1}{2}(-i\sqrt[3]{2}\sqrt{3}-\sqrt[3]{2}),\sqrt[3]{2}\right\rbrace$, no obstante, los estudiantes no encontraron este comando, por lo tanto no lo utilizaron.
    \item \textbf{Calcme: } este es un \textit{artefacto digital} que tiene una versión web, al usar el comando Resolver($z^3=2$) entrega sólo la solución real: ${{z=\sqrt[3]{2}}}$. En cambio si se utiliza el comando Resolver($z^3=2,\mathbb{C})$, es decir, especificando el conjunto sobre el que se calculan las soluciones entrega todas las soluciones complejas: 	$\left\lbrace \{z=\sqrt[3]{2} \}, \left\lbrace z=-\dfrac{\sqrt[3]{2}}{2}+\dfrac{\sqrt[6]{108} \cdot i}{2}\right\rbrace , \left\lbrace z=-\dfrac{\sqrt[3]{2}}{2}-\dfrac{\sqrt[6]{108} \cdot i}{2}  \right\rbrace  \right\rbrace $. Los estudiantes no utilizaron este software, pero se utilizó para programar las tareas en Moodle/Wiris, por lo tanto, estas soluciones aparecen en el feedback.
\end{itemize}

\begin{table}[h]
\caption{Resumen de la forma en que entregan la solución de los artefactos digitales. El * está marcado porque las soluciones complejas que entrega Wolfram Alpha son erróneas, puesto que son otras formas de escribir una la solución real}
\label{Tabla:2_soluciones_de_los_artefactos}
\begin{center}
\scalebox{0.7}{\begin{tabular}{p{2.4cm} p{2.5cm} p{2.5cm} p{2.5cm} p{2.5cm} p{2.5cm}}
\textbf{Artefacto / Soluciones}         & \textbf{Solo una solución real} & \textbf{Una   solución real y dos complejas} & \textbf{Soluciones   exactas} & \textbf{Soluciones   aproximadas} & \textbf{Represent.   gráfica de las soluciones} \\ \midrule
GeoGebra      & Sí                     & Sí                                  &                         & Sí                       & Sí                                         \\ \midrule
Symbolab      &                        & Sí                                  &                         &                          &                                            \\ \midrule
Wolfram Alpha &                        & Sí                                  & Sí*                     & Sí                       & Sí                                         \\ \midrule
Calcme        & Sí                     & Sí                                  & Sí                      &                          &                                            \\ \midrule
Photomath     & Sí                     &                                     &                         &                          & \\ \bottomrule
                                          
\end{tabular}}
\end{center}
\end{table}

Al comparar los \textit{artefactos digitales} se observa que cada uno  entrega información distinta, utilizan distintos registros (gráfico y algebraico) y en el registro algebraico entregan soluciones exactas (uno de ellos entrega información errónea) y otras aproximadas. Esto se resume en la tabla \ref{Tabla:2_soluciones_de_los_artefactos}. Cada uno arroja distintos resultados, según sus algoritmos internos, que dan cuenta de la inteligencia histórica de siglos de desarrollo sobre los números complejos. Además, el uso de aproximaciones y el signo igual ``$=$'' muestran la relatividad de la validez epistemológica a la que se enfrentan los usuarios. Este resumen no se hace con el fin de determinar cuál es el mejor artefacto, más bien, para mostrar la variedad de símbolos a los que se pueden enfrentar los estudiantes y los significados que estos pueden dar a partir de una tarea que pida una de las soluciones específica. De hecho, si quisiéramos hacer una comparación, pareciera que Symbolab es la que entrega una mejor respuesta al resolver la ecuación $z^3=2$, sin embargo, si se les pide resolver $z^5=2$, se observa que solo entrega la solución real, por lo que su validez epistemológica sigue siendo relativa.
% Please add the following required packages to your document preamble:
% \usepackage{booktabs}

\section{Resultados}

\begin{figure}[h]
    \centering
    \includegraphics[height=6.5cm]{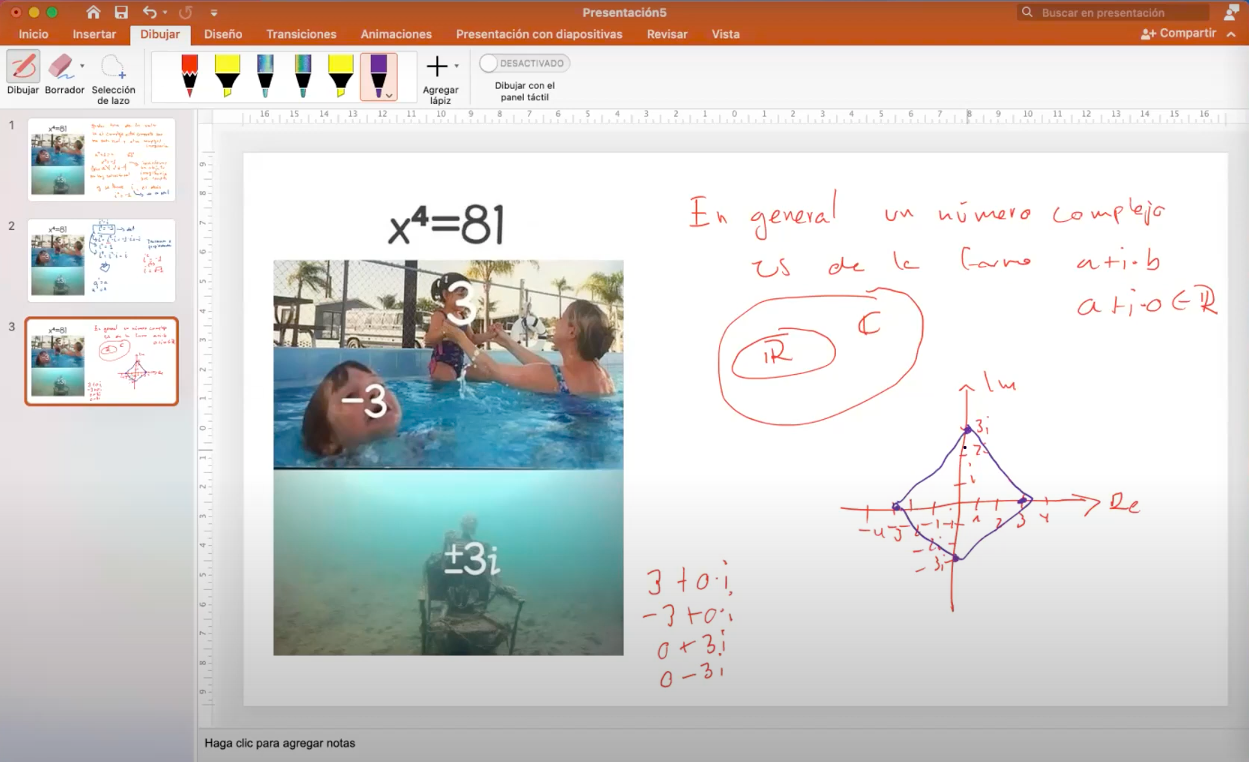}
    \caption{Meme para realizar el diagnóstico a los estudiantes sobre sus conocimientos sobre números complejos junto a la lluvia de ideas a partir de los conceptos que recordaban los estudiantes.}
    \label{fig:6_diagnostico}
\end{figure}

A continuación se muestra el \textbf{diagnóstico} que hizo el profesor, el cual consistió en preguntar a los estudiantes qué entendían de la imagen que aparece en la parte izquierda en la diapositiva de la figura \ref{fig:6_diagnostico}, con el objeto de dialogar sobre las nociones matemáticas que ellos recordaban. A medida que los estudiantes nombraban conceptos, el profesor los fue escribiendo en la diapositiva (parte derecha de la figura \ref{fig:6_diagnostico}). Una vez que hablan de la representación geométrica de los números complejos, el profesor ubica las soluciones en un plano cartesiano y pregunta ¿qué figura se forma y cómo se puede justificar que es la figura indicada? (ver figura \ref{fig:6_diagnostico}), algunos dicen que es un rombo y otros un cuadrado y se discuten algunas condiciones para justificar cada una de estas afirmaciones.

Este diagnóstico motivó a los estudiantes a movilizar sus conocimientos previos, ir conectando algunas nociones, pero también, evidenció que el \textit{referencial teórico} con el que contaban los estudiantes era débil respecto a este sistema numérico. A continuación presentamos el análisis de las tareas 1, 2 y 3, dónde hemos identificado \textit{potencialidades} y \textit{dificultades}  en el tranajo matemático orientado por artefactos digitales.

%El diagnóstico se realizó mostrando en pantalla la imagen que aparece en la parte izquierda en la diapositiva de la figura \ref{fig:6_diagnostico}. Se les preguntó a los estudiantes qué entendían de ella.  Entregaron respuestas sobre algunas nociones superficiales y aisladas sobre números complejos. A medida que los estudiantes nombraban conceptos, el profesor los fue escribiendo en la diapositiva (parte derecha de la figura \ref{fig:6_diagnostico}). Esto permitió ir conectando algunos conceptos y a otros estudiantes recordar elementos, pero también, evidenció que el \textit{referencial teórico} con el que contaban los estudiantes era débil respecto a este sistema numérico. 

%Este diagnóstico motivó a los estudiantes a movilizar sus conocimientos previos sobre estas nociones, y da el punto de partida para la puesta en práctica de las tareas que siguen durante las sesiones de clases, su análisis se presenta a continuación. 

\subsection{Potencial del trabajo matemático orientado por los \textit{artefactos digitales} }

\subsubsection{Las limitaciones instrumentales que activan la \textit{génesis discursiva}}

\begin{figure}[h]
    \centering
    \includegraphics[width=0.9\textwidth]{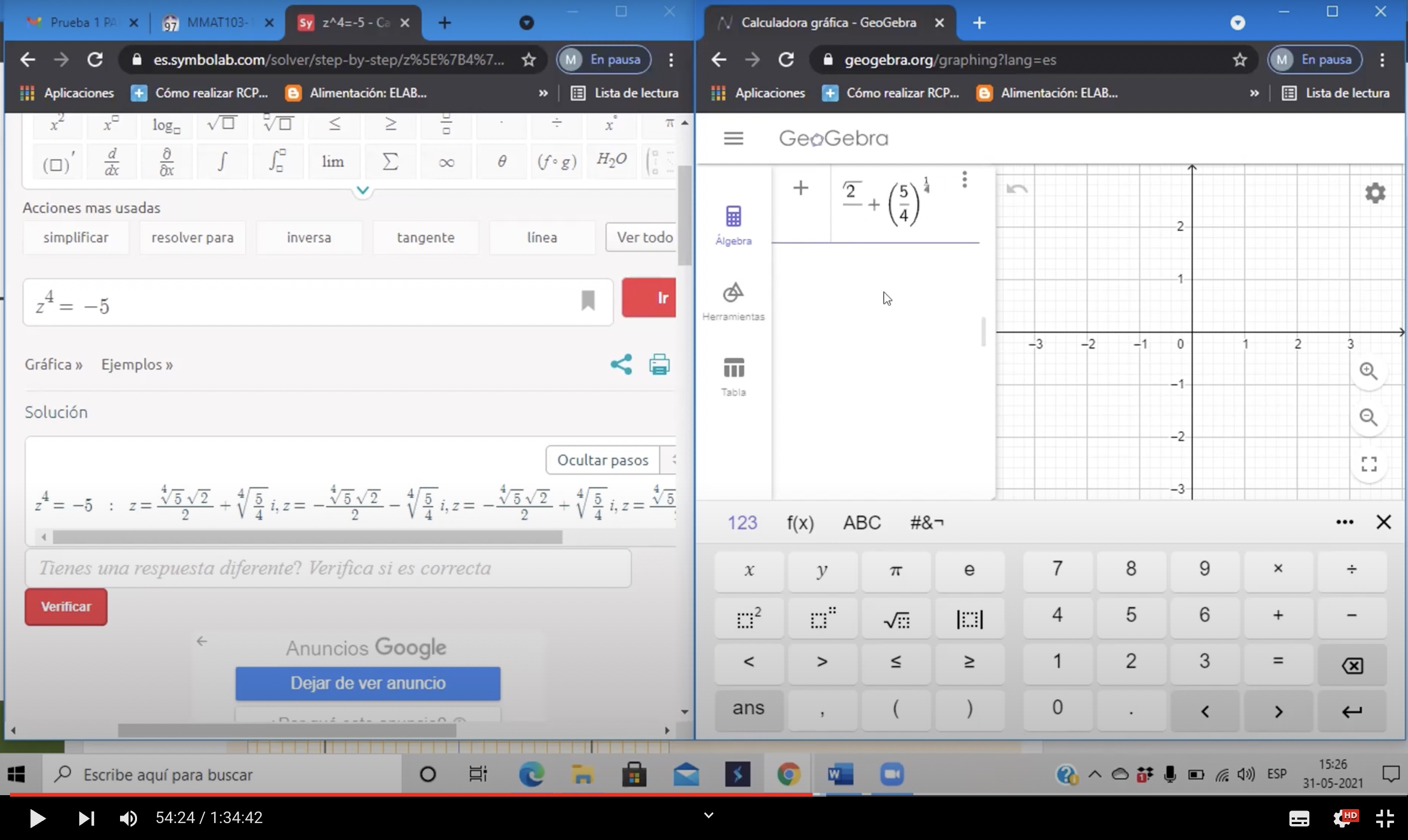}
    \caption{Captura de pantalla del trabajo de Sofia ingresando una potencia como raíz}
    \label{fig:9_Sofia_ingresando_raiz}
\end{figure}

De la tarea dos, obtenidas las soluciones de la ecuación $z^4=-5$, seleccionamos el episodio de Sofía quien señala que no entiende (minuto 46:26 tabla \ref{tab:Potencialidad_dialogo_Sofia_ingresa_soluciones}) la tarea y en particular el significado de las soluciones dadas por Symbolab (el \textit{artefacto digital} usado). Posteriormente el profesor le pide que abra la calculadora de GeoGebra y que en paralelo comparta pantalla con Symbolab, para mostrar lo que el \textit{artefacto digital} da como soluciones para la ecuación (ver figura \ref{fig:9_Sofia_ingresando_raiz} y tabla \ref{tab:Potencialidad_dialogo_Sofia_ingresa_soluciones}). 

\begin{table}[]
\caption{Transcripción del diálogo donde Sofia ingresa la solución de Symbolab en Geogebra}
\label{tab:Potencialidad_dialogo_Sofia_ingresa_soluciones}
\scalebox{0.6}{\begin{tabular}{p{20cm}l@{}}
\toprule
Sofia: 46:26 Profe, ¿sabe qué? Sinceramente, no sé cómo hacerlo.   Como que me explota el cerebro. O sea, estoy ocupando Symbolab, pero… o sea,   como que no… no entiendo en sí el ejercicio. No, no logro entender que por   qué da eso.                                                                                                                                                        \\ \midrule
Sofia: 47:03 Ya, mire. Eso es lo que sale [mostrando lo que le entrega Symbolab], pero no sé qué es lo   que me sirve y qué es lo que no, y por qué da eso.                                                                                                                                                                                                                                                                              \\ \midrule
Profesor: 47:17 Mira, ¿puedes abrir otra pestaña donde abras la   calculadora de Geogebra? Por favor. Y después… no sé si sabes separar las pestañas   para que quede una al lado de la otra.                                                                                                                                                                                                                                                                                           \\ \midrule
Sofia: 48:37 ¿Y cómo? ¿Cómo comparto dos?                                                                                                                                                                                                                                                                                                                                                        \\ \midrule
Profesor: 48:40 No, tienes que compartir el escritorio. Una de   las opciones cuando compartes es compartir todo el escritorio                                                                                                                                                                                                                                                                     \\ \midrule
Profesor: 49:07 Ahí sí. Entonces vamos… lo que vamos a hacer es   graficar todas esas soluciones que tienes ahí. Porque lees y lo que no se entiende es que significan todas eso símbolos raros que están   ahí, las raíces y todo eso?.                                                                     \\ \midrule
Profesor: 49:31 Entonces, en la primera dice raíz cuarta de 4   por raíz de 2, partido por 2, más la raíz cuarta de 5/4 por ¿Cierto?. Ingrésala ahí al Geogebra                                                                                                                       \\ \midrule
Sofia: 49:54 Ya, copiar y pegar no más. Oh, no se puede.      \\ \midrule
Sofia: 50:09 Ya, pero no sé cómo ingresarlo, profe.                                                                            \\ \midrule
Sofia: 50:14 Tengo que poner $z=$. No sé cómo se hace           \\ \midrule
Profesor: 50:24 Tranquila, si para eso estamos en este tipo de cursos. Para aprender a usar este tipo de software. A ver ¿cómo se puede ingresar la raíz cuarta si es que ahí no te aparece en el menú? ¿Hay alguna   otra manera de escribir raíces? ¿Tú sabes qué significa la raíz cuarta? ¿Cómo   se puede escribir de otra forma? ¿O alguien puede ayudar si es que Sofia no   lo sabe?  \\ \midrule

Cristian: 50:51 Cómo potencias.     \\ \midrule

Profesor: 50:52 Con potencias,¿la raíz cuarta cómo se   escribiría como potencias? Por ejemplo, la raíz cuarta de 5, ¿cómo se puede escribir como   potencias?                                                                                                     \\ \midrule
Juan: 50:56 Es 5 elevado a 1/4.                                                                                                                                                                 \\ \midrule
Profesor: 50:59 Ya, entonces puedes escribirlo así Sofia.                                                                    \\ \midrule
Sofia: 51:01 Oh, quedé máquina. Ya, 5 elevado a ¿y cómo   escribo?                                                                                                                      \\ \midrule
Sofia: 51:13 Uh, que loco. Ya. ¿Y después cómo lo escribo   partido en 2?                                                   \\ \midrule
Profesor: 51:23 Ahora tienes que multiplicarlo por raíz de 2.                                                                  \\ \midrule
Profesor: 51:49 Por raíz de 2 dice ahí. No dice por 2. Ojo.                                                                    \\ \midrule
Sofia: 51:55 Ah, verdad, verdad.                                \\ \midrule
Sofia: 52:19 Oh, que magia, profe. Estoy impactada.                                                                            \\ \midrule
Profesor: 52:28 ¿Y ahora como escribiría en la raíz cuarta de   5/4?                                                                                                                          \\ \midrule
Sofia: 52:41 No sé, profe, porque. No sé escribir esa raíz de 5/4. O sea, la   cuarta raíz de 5/4.                                                                      \\ \midrule
Profesor: 53:21 Ya, ¿alguien podría darle una recomendación a   Sofia, cómo hacerlo? ¿Cómo escribir la raíz cuarta de 5/4 con lo que   tenemos ahí a la vista?                            \\ \midrule
Sofia: 53:37 Sería 5/4 elevado a 1/4.                                                                                            \\ \midrule
Profesor: 53:41 Muy bien, Sofia, máquina.                                                                                    \\ \midrule
Sofia: 53:40 Soy una crack. Ay, yo pensé que influía algo el   5/4. Ah, pero no.                                               \\ \midrule
Profesor: 54:25 Ya, y después eso lo multiplicas por i. \\ \midrule        \\ Sofia: 54:00 Listo profe                                                                                                       \\ \bottomrule
\end{tabular}
}
\end{table}

 Posterioremente, Sofía ingresa las soluciones obtenidas en Symblab a GeoGebra para graficarlas, en este momento hay dificultades para utilizar el artefacto en la escritura de raíces (minuto 50:14 tabla \ref{tab:Potencialidad_dialogo_Sofia_ingresa_soluciones}, ya que GeoGebra no tiene un comando que permita la escritura directa de la raíz con un índice distinto a 2. Sofía no recuerda de qué otra forma ingresar la solución $z=\sqrt[4]{5}\sqrt{2}+\sqrt[4]{\dfrac{5}{4}}i$, sus compañeros y el profesor retroalimentan a Sofia, y le señalan que debe utilizar la expresión usando notación de potencias para ingresar las raíces en GeoGebra. A partir de esta retroalimentación, Sofía logra ingresar las raíces como potencias activando no sólo una cuestión notacional, si no que conocimientos asociados a nociones de potencias y radicales.

En términos teóricos podemos observar que se activa la \textit{génesis instrumental} en el caso de Sofía mediante el uso de Symbolab, los estudiantes observan las soluciones dadas por el \textit{artefacto digital}, pero no existe todavía una \textit{visualización} de estos símbolos, ya que hasta ese momento no se articula el signo con su significado, o dicho de otra manera existe una visualización icónica del signo. 

Mediante la retroalimentación del profesor y los estudiantes, se logra superar dificultades \textit{instrumentales}, que permite activar la \textit{génesis discursiva}, ya que la estudiante comprende que es necesario usar propiedades de los radicales para ingresar las raíces con índice distinto a dos en GeoGebra, observamos entonces activación del plano \textit{[Ins-Dis]}, para ingresar las soluciones dadas mediante propiedades de potencias.

%Para intentar dar sentido a esas soluciones, el profesor le pide que trabaje en paralelo con ambos artefactos (Symbolab y GeoGebra). En esta parte emergen dificultades \textit{instrumentales} para ingresar las soluciones que da Symbolab en GeoGebra pero, se logran superar gracias a la retroalimentación de los compañeros y del profesor. En este momento,  se activa la \textit{génesis discursiva}, ya que la estudiante comprende que es necesario usar propiedades de potencias para ingresar las raíces con índice distinto a dos en GeoGebra, observamos entonces activación del plano \textit{[Ins-Dis]}, para ingresar las soluciones dadas mediante propiedades de potencias.

\subsubsection{La articulación de tres \textit{artefactos digitales} ayuda a la visualización}

\begin{figure}[h!]
    \centering
    \includegraphics[width=0.8\textwidth]{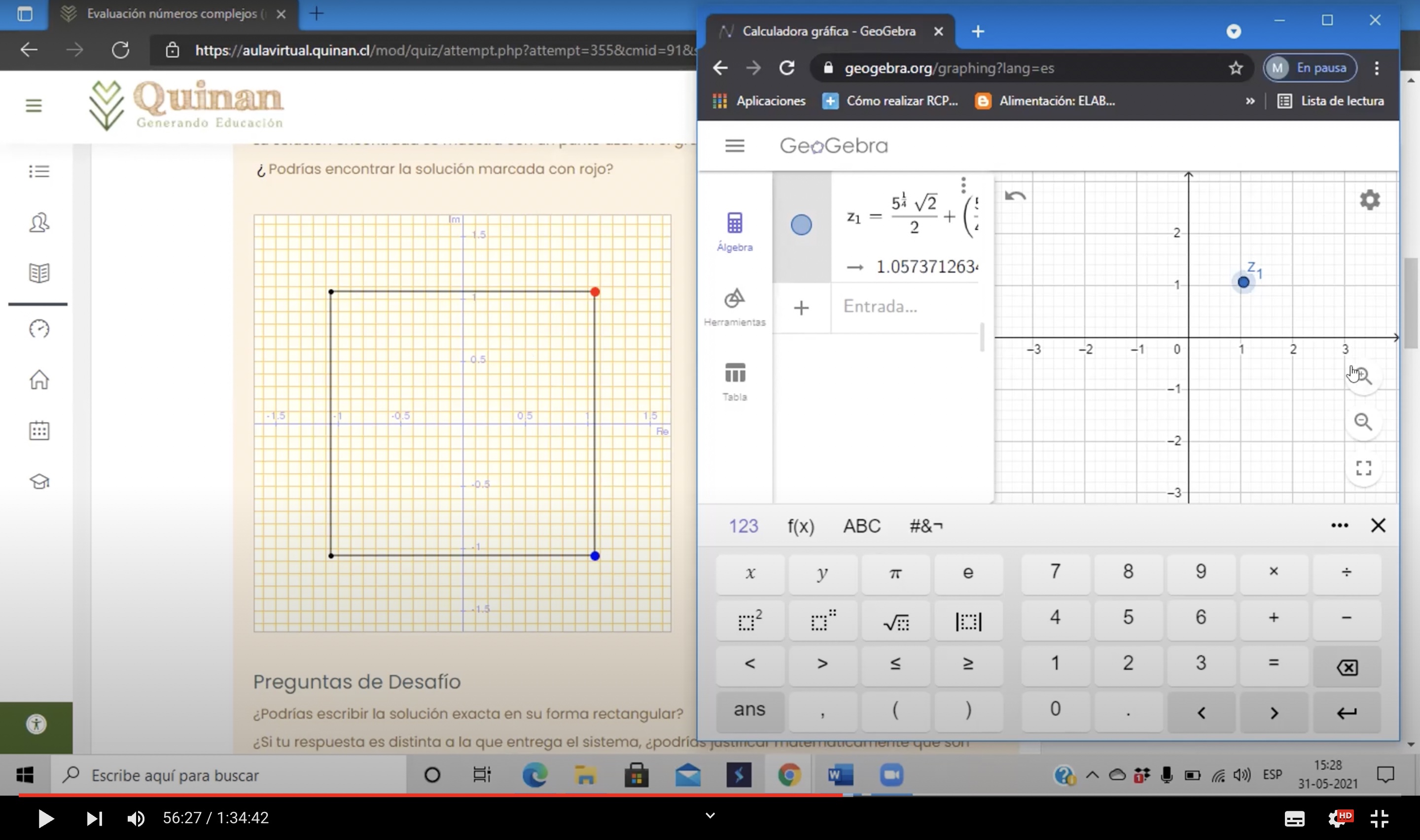}
    \caption{Captura de pantalla del trabajo de Sofia donde compara la información que entrega la tarea en Moodle/Wiris y el punto que graficó en GeoGebra.}
    \label{fig:10_Sofia_articulando_dos_registros}
\end{figure}

Nuevamente destacamos el episodio de Sofía, quien  luego de haber ingresado las soluciones en GeoGebra, éste las representa en el plano cartesiano (dibuja un punto), la estudiante comparte la ventana de Moodle/Wiris (la tarea) y GeoGebra (la solución) (ver figura \ref{fig:10_Sofia_articulando_dos_registros} y \ref{tab:potencialidades_Sofia_comprende}), en este momento visualiza las raíces de la ecuación en su representación gráfica logrando comprender que el punto rojo pedido en la tarea es la solución correcta. Nos parece que este episodio es altamente significativo, ya que previamente, las soluciones entregadas por Symbolab no tenían sentido para ella, puesto que no existía la relación de la ecuación dada en el registro algebraico y el cuadrado presentado en la tarea.\\

%\scalebox{0.7}{\begin{tabular}{p{18cm}l@{}}
% Please add the following required packages to your document preamble:
% \usepackage{booktabs}
% Please add the following required packages to your document preamble:
% \usepackage{booktabs}

\begin{table}[]
\caption{Continuación Trascripción del diálogo de Sofia comparando la información que entrega la tarea en Moodle/Wiris y el punto que graficó en GeoGebra. }
\label{tab:potencialidades_Sofia_comprende}
\scalebox{0.6}{\begin{tabular}{p{20cm}l@{}}
\toprule
Sofia: 46:26 Profe, ¿sabe qué? Sinceramente, no sé cómo hacerlo.   Como que me explota el cerebro. O sea, estoy ocupando Symbolab, pero… o sea,   como que no… no entiendo en sí el ejercicio. No, no logro entender que por   qué da eso.                                                                                                                                                        \\ \midrule
...                                                                                                                                                                     \\ \midrule
Sofia: 55:10 Entonces, ese es un punto que está… que sería…   sería el de acá, ¿o no? No.                                     \\ \midrule
Sofia: 55:32 Ya me dí cuenta.                                                                                                  \\ \midrule
Sofia: 55:38 O no sé si… a ver, a ver déjeme seguir analizando.   Es el punto rojo, ¿o no? Este de acá.                          \\ \midrule
Sofia: 55:56 Porque aquí están los reales y los imaginarios   están acá y aquí está el 1 en positivo y aquí también. No sé cómo fundamentar   mi respuesta, profe, la verdad.              \\ \midrule
Sofia: 56:18 O sea, los relaciono más que nada.                                                                                \\ \midrule
Sofia: 56:26 ¿Aquí?                                                                                                            \\ \midrule
Profesor: 56:29  Para ver…   el gráfico. Ahí sale un signo +, una lupita. ¿Qué le pasa a lupita?                              \\ \midrule
Profesor: 56:54 ¿Ahora estás entendiendo lo que estás buscando?   ¿no?                                                                                                                         \\ \midrule
Sofia: 56:57 Sí, profe, ahora me queda completamente claro. ¿Y   qué es lo que te pide acá? Ahora hay que entender lo que te pide.                                                           \\ \midrule
Sofia: 58:44 Ah, es que esta es una suma, profe.                                                                                                                                              \\ \midrule
Sofia: 59:55 Entendí cómo pasar estas cosas. Las raíces, o sea,   la cuarta raíz de 5 a potencia.                                                                                               \\ \midrule
Sofia: 60:12 Pucha, Geogebra. Lo adoro. Tengo que aprender a   usarlo sí. Y que la diferencia en la respuesta que había puesto era que el   otro era una resta y esta era una suma. Porque el cuadrante de acá es   positivo y el otro es negativo, por eso encontré el punto negativo de acá. No   sé si está bien lo que dije.                                                                  \\ \bottomrule
\end{tabular}
}
\end{table}

%En este momento Sofia logra dar significado a las soluciones dadas en su forma algebraica, rectangular y exacta, ya que al ingresarlas al GeoGebra éste artefacto las representa en el plano cartesiano (dibuja un punto), y finalmente ella logra comprender el punto rojo pedido en la tarea, lo cual previamente no tenía sentido para ella, ya que no comprendía la relación de la ecuación dada en el registro algebraico con el cuadrado que se presentaba en el plano cartesiano. 

%Esto no solo ayudó a Sofia, también a otros compañeros, como Robert, quien indica lo siguiente: ``EH: 61:30 profe a mí me ayudó lo que dijo Sofia en un momento son los signos de cada parte de un binomio representa al cuadrante donde están entonces la primera parte del binomio me está diciendo que está en los positivos de los reales y la segunda parte en los negativos de los reales, en este caso la corrección de esto sería que la resta debió haber sido una suma. '' lo cual da muestra del potencial de la discusión abierta en  clases.

En términos teóricos, en este episodio se activa el plano \textit{[Sem-Ins]}, ya que al utilizar el GeoGebra para \textit{visualizar} las soluciones de $z^4=-5$ obtenidas previamente por Symbolab, se logra dar sentido a las raíces de la ecuación mediante la articulación de los registros algebraico y gráfico.  Lo anterior, le ha permitido comprender las respuestas entregadas por Symbolab, esto revela que la activación de solamente la \textit{génesis} \textit{instrumental}, no es suficiente para entender los objetos matemáticos involucrados.

%Luego, utiliza \textit{artefacto digital} GeoGebra para representar una solución sobre el plano complejo, coincidiendo con la solución pedida en la tarea (figura \ref{fig:10_Sofia_articulando_dos_registros}, lado izquierdo). En este momento, ella logra dar sentido a las soluciones de la ecuación mediante la articulación de los registros algebraico y gráfico. Lo anterior, le ha permitido comprender las respuestas entregadas por Symbolab, esto revela que la activación de solamente la \textit{génesis} (\textit{instrumental}), no es suficiente para entender los objetos matemáticos involucrados, por lo que es necesario que haya una activación del plano \textit{[Sem-Ins]} para que la estudiante comprenda que cada una de las soluciones es un punto en plano cartesiano. Esto, además le permitió corregir su error previo, ya que ella había ingresado una solución, pero no era la pedida en la tarea, ahora Sofia logra determinar la solución correcta mediante la activación del plano vertical \textit{[Sem-Ins]}.

\subsubsection{El \textit{artefacto digital} guía el trabajo \textit{discursivo}}
\begin{figure}[h!]
    \centering
    \includegraphics[width=0.99\textwidth]{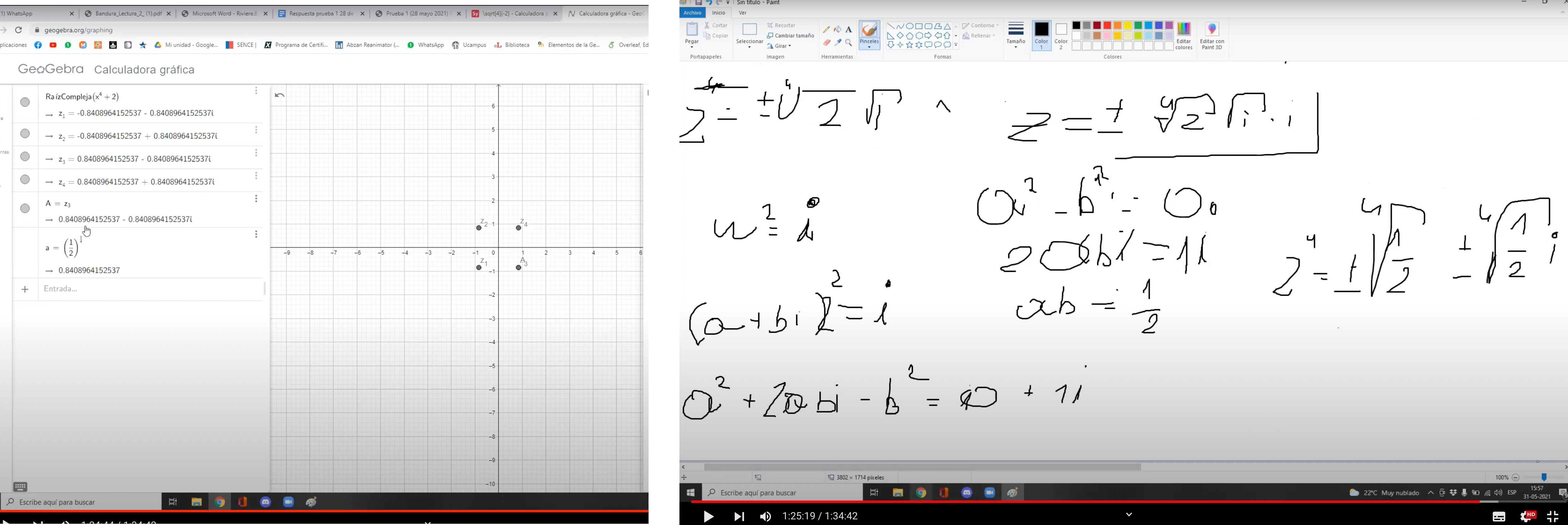}
    \caption{Capturas de pantalla de la explicación de Juan}
    \label{fig:11_Argumentacion_de-Jose_Ignacio}
\end{figure}

Otro caso seleccionado es el episodio de Juan en la tarea dos, quien comparte pantalla (y explica en una pizarra de Paint\footnote{Software de dibujo en el sistema operativo Microsoft}), su trabajo matemático al resolver la ecuación $z^4=-2$. Explica  que al utilizar en GeoGebra el comando ``RaízCompleja'' es necesario reescribir la ecuación de forma polinómica igualando a cero: $RaizCompleja(z^4+2)$ en vez de $RaizCompleja(z^4=-2)$. GeoGebra muestra las cuatro soluciones aproximadas (por ejemplo: $z_1-0.8408...-0.8408i$), y por medio de una estrategia de ensayo y error evalúa las soluciones en la ecuación para verificar si la satisfacen (ver figura \ref{fig:11_Argumentacion_de-Jose_Ignacio} parte izquierda).

Posterior a eso, Juan cambia de estrategia y resuelve la ecuación de forma algebraica. Esta  estrategia consiste en aplicar reiteradamente raíz cuadrada a toda la ecuación, (Juan: 76:28, tabla \ref{tab:5_transcripcion_Jose_Ignacio}) obteniendo $z=\pm \sqrt[4]{2}\sqrt{i}\cdot i$. El profesor  interviene preguntando por el significado de  $\sqrt{i}$ (Profesor: 80:36, tabla \ref{tab:5_transcripcion_Jose_Ignacio}), lo cual da paso al estudiante a argumentar cómo se puede escribir raíz de $i$ de forma canónica ($a+bi$) (Juan: 80:39-84:14, tabla \ref{tab:5_transcripcion_Jose_Ignacio}), utiliza un sistema de ecuaciones no lineal ($a^2-b^2=0$ y $2ab=1$), que es resuelto por Wolfram Alpha por recomendación del profesor. Finalmente, Juan retoma las soluciones aproximadas entregadas por el GeoGebra, y las compara con el trabajo algebraico que realizó para argumentar cómo se llega estos valores. (Juan: 88:57, tabla \ref{tab:5_transcripcion_Jose_Ignacio}).   

% Please add the following required packages to your document preamble:
% \usepackage{booktabs}

\begin{table}[]
\caption{Transcripción de la explicación algebraica de Juan}
\label{tab:5_transcripcion_Jose_Ignacio}
\begin{center}
\scalebox{0.7}{\begin{tabular}{p{18cm}l@{}}
\toprule
Juan: 76:28 al ingresar la ecuación $z^4+2$ a   GeoGebra, éste nos muestra los valores en forma de número como para darme una   idea...después lo abordamos desde una perspectiva algebraica. Lo hice con un   amigo, lo que hicimos fue aplicamos una raíz   cuadrada para que nos quedara una expresión de este estilo: más menos raíz de   menos dos $z^2=\pm \sqrt{-2}$. Eso es lo que nos quedaría al aplicar una raíz cuadrada. Entonces, bueno,   aquí básicamente estos son dos valores, el positivo y el negativo, $z$ es igual   a más raíz cuadrada de menos dos $z=\sqrt{-2}$, y lo mismo en el otro lado, menos raíz cuadrada de menos dos $z=-\sqrt{-2}$. Esto   lo podemos expresar de forma compleja, esto   es lo mismo que decir dos por menos uno, separándolo queda $z$ cuadrado es   igual a raíz de dos i ${[}z^2=\sqrt{2i}{]}$ y $z$ cuadrado es igual a menos raíz de $2i$   ${[}z^2=-\sqrt{2i}{]}$, quedando dos expresiones cuadradas, luego a estas dos le   aplica raíz cuadrada y queda como, $z$ es igual a más menos raíz cuarta de dos   por raíz cuadrada de $i$ ${[}z=\pm \sqrt[4]{2}\sqrt{i}{]}$, y ahí aparecen los   cuatro valores. \\      
\midrule
Profesor:   80:36 ¿Y la raíz de i que significa?                                                                           \\\midrule
Juan: 80:39 bueno ese fue el siguiente desafío,   determinar cuánto vale la raíz de i, bueno básicamente es raíz de menos uno $i=\sqrt{-1}$, ayudándonos de un valor auxiliar, para decir por   ejemplo, ``$w^2=i$'' lo expresamos en forma de polinomio porque esto es un   número complejo, y luego decimos esto vale $(a+bi)^2=i$.   Desarrollamos quedando $a^2+2abi-b^2=0+i$, luego iguala la parte real a   cero y la parte imaginaria a $i$. Entonces como entendí yo de lo que hizo mi   amigo, a cuadrado menos b cuadrado  tenía que ser igual a cero ${[}a^2-b^2=0{]}$,   y dos por abi igual a 1 por $i$ ${[}2abi=1i{]}$, bueno el dos pasa al   otro lado y nos quedaría como ab porque las i se eliminan   ¿no?                                                      \\ \midrule
Profesor:   83:16 Claro, las $i$ se eliminan, por definición                                                   \\ \midrule
Juan: 84:14 ya bueno, y a ver, aquí recordando con   GeoGebra habíamos llegado a que ese valor, que es el 0.8408. Era lo mismo   que decir raíz cuarta de un medio $ \left[ \sqrt[4]{\dfrac{1}{2}} \right]$, 
entonces ya más o menos sabíamos que teníamos que   llegar a una expresión del estilo $z$ es igual a más menos raíz cuarta de un medio y   más menos lo mismo, pero con $i$: $ \left[z=\pm \sqrt[4]{\dfrac{1}{2}} \pm   \sqrt[4]{\dfrac{1}{2}}i\right]$, 
sabíamos que teníamos que llegar a eso más o menos   por GeoGebra, llegamos a esa conclusión. Si soy sincero, no me acuerdo cómo   fue que logramos resolver esto \\ \midrule
Profesor: 85:31 Es que ahí aparece un sistema. Tienes $a^2–b^2=0$ y $2a=1$. Entonces es cosa de resolver ese sistema que aparece ahí nomás. De hecho, abre Wolfram, está en inglés si, y coloca solve espacio y las ecuaciones y dale enter, y de ahí aparecen las soluciones, para calcular la raíz cuadrada del i $[\sqrt{i}]$ \\ \midrule

Juan: 88:57 claro!, Y con eso ya con eso, y después sólo reemplazo $i$ por la respuesta, que en este caso era una de las respuestas del sistema de ecuaciones, que es justamente el resultado al que llega GeoGebra que es menos raíz cuarta de un medio $\left[- \sqrt[4]{\dfrac{1}{2}}\right] $ o más raíz cuarta de $\left[ \sqrt[4]{\dfrac{1}{2}}\right] $ y así sucesivamente, así llegué al resultado 

    \\ \bottomrule
\end{tabular}}
\end{center}
\end{table}

 En términos teóricos se activa el plano \textit{[Sem-Ins]} porque son \textit{visualizadas} las soluciones en GeoGebra, esto ayuda a dirigir el trabajo algebraico, dando paso a la activación de la \textit{génesis discursiva} que le permite justificar los resultados obtenidos por el \textit{artefacto digital}. El estudiante realiza todo el trabajo \textit{discursivo} que no realizan los \textit{artefactos digitales}, que a su vez, es controlado mediante el plano \textit{[Sem-Ins]} a partir del uso del \textit{artefacto digital} (GeoGebra), ya que verifica y valida lo que se va obteniendo algebraicamente apoyado por las  soluciones de GeoGebra.

\subsection{Dificultades en el trabajo matemático con \textit{artefactos digitales}}

\subsubsection{Dificultades con el sistema de validación del CAA.}
\begin{figure}[h]
    \centering
    \includegraphics[height=6.5cm]{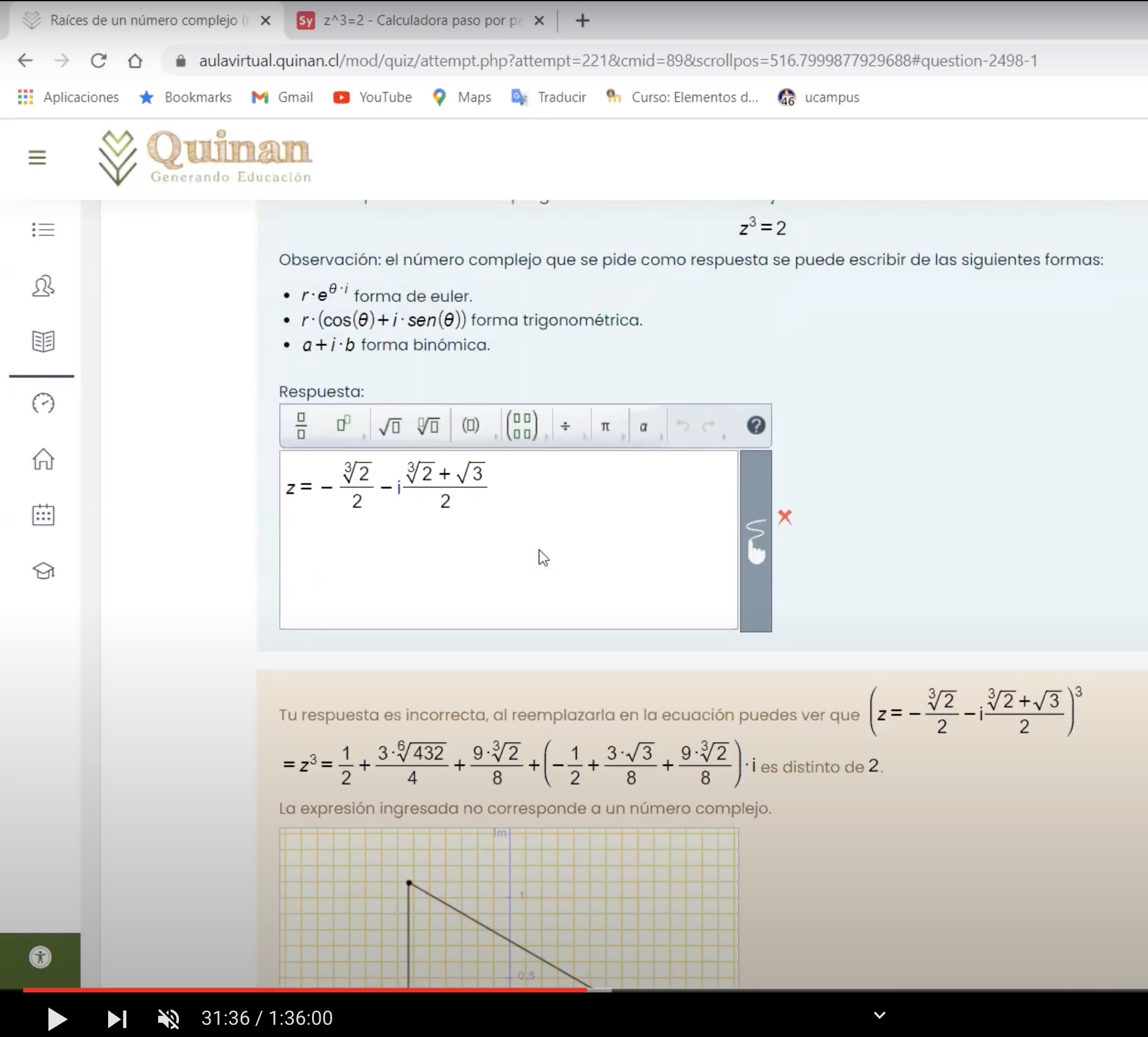}
    \caption{Feedback entregado por Moodle/Wiris.}
    \label{fig:7_dificultad_ingresar_respuestas}
\end{figure}

De la tarea uno, seleccionamos la producción de Paula quien utiliza Symbolab, este artefacto le arroja tres soluciones algebraicas, rectangulares y exactas. A pesar que Paula identifica correctamente la solución de la tarea (vértice del triángulo dado), ella tiene una dificultad para ingresar la solución en Moodle/Wiris, ya que incluye ``$z=$'', pero el Moodle/Wiris admite sólo el lado derecho de la igualdad $z=a+ib$, luego el sistema se la considera incorrecta (ver figura \ref{fig:7_dificultad_ingresar_respuestas}).  

\begin{table}[h]
\caption{Transcripción de diálogo de Paula con el profesor }
\label{tab:3_Dialogo_Paula_Profesor}
\begin{center}
\scalebox{0.7}{
\begin{tabular}{p{18cm}}
\toprule
Paula: 31:20 saqué de acá la respuesta como los dos son negativos van quedar en el cuadrante de abajo. No sé el nombre de los cuadrantes, perdón y puse aquí la respuesta y me la tomó como mala porque estos no estaban como sumados. \\ \midrule

Profesor: 31:37 el problema es que va sin el $z=$, ingresa de nuevo, coloca arriba finalizar revisión y has un nuevo intento.                                                                                                                                                            \\  \bottomrule                                                                      
\end{tabular}
}
\end{center}
\end{table}

Podemos observar que el feedback no es suficiente para superar la dificultad de escritura de la solución, y es el profesor quien tiene que indicarle porqué Moodle/Wiris considera incorrecta su respuesta. Este caso es relevante pues permite tener una retroalimentación al resto de los estudiantes, para que tengan conciencia sobre cómo deben ser ingresadas las soluciones en Moodle/Wiris.
\\

En términos teóricos, en este episiodio los estudiantes activan el plano \textit{[Sem-Ins]}, ya que utilizan un artefacto digital (Symbolab u otro) que entrega las tres soluciones, los estudiantes  identifican la solución correcta solicitada por Moodle/Wiris estableciendo una relación entre los cuadrantes y el signo de las soluciones. Sin embargo, debido al formato específico en que se pide la respuesta, es decir, una dificultad instrumental, los estudiantes que tienen esta dificultad (como el caso de Paula) podrían haber quedado con una interpretación errónea a pesar de haber razonado de forma correcta. 

\newpage

\subsubsection{Dificultades cuando el \textit{artefacto digital} presenta un error en la respuesta dada}

\begin{figure}[h!]
    \centering
    \includegraphics[width=0.99\textwidth]{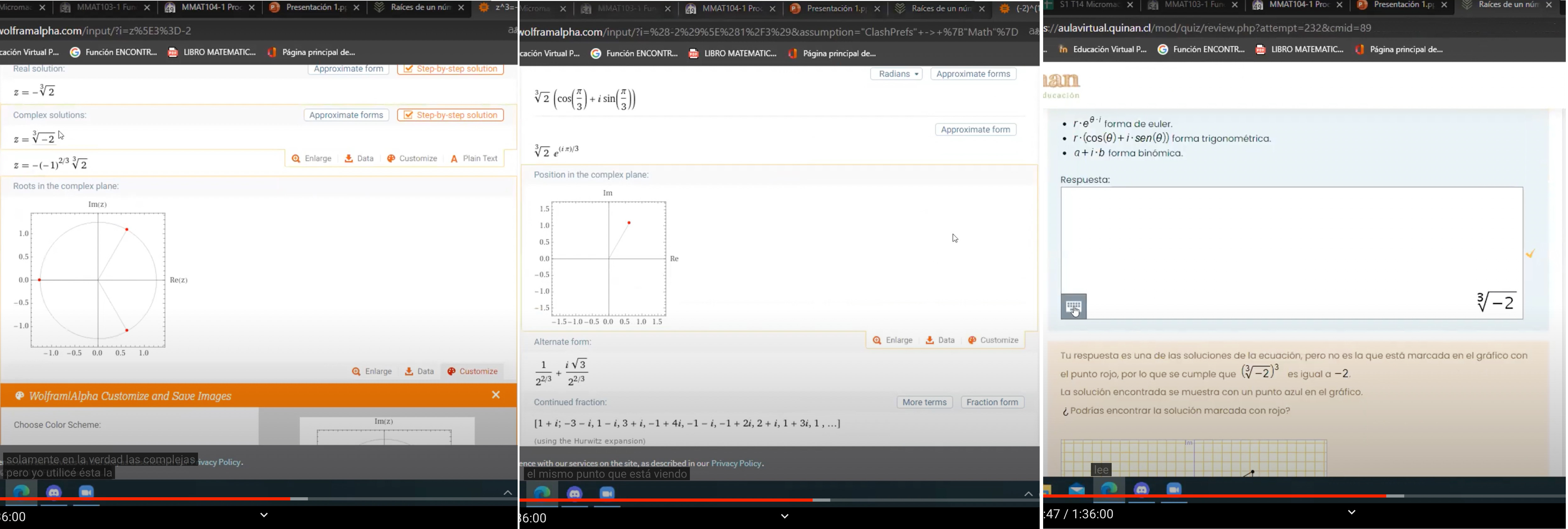}
    \caption{Captura de pantalla de distintos momentos donde Cristian explica lo que hizo en Wolfram Alpha (figura izquierda y centro) y el feedback que le entregó Moodle/Wiris al ingresar su respuesta (figura derecha)}
    \label{fig:8_explicacion_de_Cristian}
\end{figure}

De la tarea uno, seleccionamos la producción de Cristian quien utiliza Wolfram Alpha. Pero este \textit{artefacto digital} da una respuesta incorrecta, puesto que, dice: la solución real es: $z=-\sqrt[3]{2}$ y las soluciones complejas son: $z=\sqrt[3]{-2}$ y $z=(-1)^{2/3}\sqrt[3]{2}$. También, muestra las soluciones en un plano complejo (ver tabla \ref{tab:3_Discusion_con_Cristian} y figura \ref{fig:8_explicacion_de_Cristian}). Cristian elije una de estas respuestas ($z=\sqrt[3]{-2}$) tomando en cuenta la información gráfica, pero Moodle/Wiris le indica que es una de las soluciones de la ecuación, pero no la del vértice solicitado.

\begin{table}[h!]
\caption{Transcripción de discusión con Cristian sobre el uso de Wolfram Alpha, junto con pantallazos mostrados por el estudiante durante la conversación con el profesor.}
\label{tab:3_Discusion_con_Cristian}
\begin{center}
\scalebox{0.7}{
\begin{tabular}{p{18cm}}
\toprule
Cristian:   50:09 {[}…{]} escribí: $z$ elevado a 3 es igual a menos 2 {[}$z^3=-2${]} entonces están   todos los resultados acá, la solución real y las dos complejas, entonces aquí   me aparece el gráfico, y aquí es la solución real porque está en la recta de   los reales, y éstas deben ser las complejas {[}muestra sobre el software las   soluciones complejas y las diferencia de la real{]}, pero yo utilicé ésta la   que no era {[}utiliza una de las soluciones erróneas que da el Wolfram Alpha{]}    entonces me iba por acá y con que Wolfram Alpha  me lanzó el mismo punto   que está pidiendo, pero al final no era... \\ \midrule
Cristian:   51:31 en Wolfram Alpha me aparece que en el plano aparece lo que me está   solicitando y aquí al final no era nada                                                                                                                                                             \\  \midrule
Profesor: 51:34 Sube un poquito, mira lee lo que dice ahí, la retroalimentación de   la primera línea anda leyéndola con calma, lo que está en   amarillo                                                                                                                              \\  \midrule

Cristian:   52:25 [leyendo] tu respuesta es una de las soluciones de la ecuación, pero no es la que   está marcada en el gráfico con el punto rojo, por lo que se cumple que esa  solución dada es igual a menos 2 {[}lee el enunciado de la plataforma{]}                                                                                                       \\  \midrule
Cristian:   52:35 profe, pero ahí hay un detalle, ¿Porqué la elevó al cubo esa expresión, si eso no es lo que escribió?                                           \\  \midrule
Cristian   52:43 porque $z$ al cubo es igual a esto {[}refiriéndose a la ecuación   $z^3=-2${]}                                                                 \\  \midrule
Profesor: 52:46 exactamente lo que dice por lo que se cumple que la raíz cúbica dice   su respuesta es una solución de la ecuación, pero no es la que está marcada   en el gráfico y: ¿porqué es una solución de la ecuación?, porque al cubo da 2 que   es lo que se está solicitando, y la solución que tu ingresaste es la que está   marcada en azul [lee el feedback que aparece en Moodle/Wiris, ver figura \ref{fig:8_explicacion_de_Cristian}]                           \\ \bottomrule                                                                      
\end{tabular}
}
\end{center}
\end{table}

En términos teóricos, en este episiodio los estudiantes activan el plano \textit{[Sem-Ins]}, ya que utilizan Wolfram Alpha e identifican cuál es la solución que le piden. No obstante, el error del artefacto digital bloquea su trabajo matemático.  El profesor tampoco se da cuenta del error, porque asume que la solución que ingresó es la real y no una de las complejas. 

Destacamos hasta qué punto pueden ser problemáticos los resultados obtenidos en un \textit{artefacto digital} específico, más aún cuando este artefacto digital es usado de forma profesional en ámbitos de la ingeniería y ciencias, en el cual no existe gran cuestionamientos sobre la validez de sus resultados.

\subsubsection{Dificultades \textit{discursivas} al interpretar un feedback del Moodle/Wiris }

%Por un lado el Symbolab entrega como solución una expresión que contiene una raíz cuadrada multiplicada por una raíz cúbica ($\sqrt[3]{2} \cdot \sqrt{3}$), en cambio, Moodle/Wiris entrega la expresión con una raíz sexta ($\sqrt[6]{108}$), y los estudiantes no logran darse cuenta que  no logran justificar el feedback que ambas expresiones son iguales. 

%A medida que transcurre el tiempo de la clase, se van resolviendo algunos problemas \textit{instrumentales} con los distintos \textit{artefactos digitales} usados por los estudiantes. Una vez que comienzan a responder correctamente aparecen nuevas interrogantes. Symbolab entrega como solución una expresión que contiene una raíz cuadrada multiplicada por una raíz cúbica ($\sqrt[3]{2} \cdot \sqrt{3}$), en cambio, Moodle/Wiris entrega la expresión con una raíz sexta ($\sqrt[6]{108}$). Una de las preguntas que plantea el feedback es si pueden justificar que ambas expresiones son equivalentes. 

%Paula indica que sí se puede justificar utilizando  propiedades de potencias. Sin embargo, al avanzar 

De la tarea uno, seleccionamos el episodio del diálogo entre Paula, Cristopher y Juan que genera una discusión en el trabajo matemático por una dificultad de factorización que se muestra en Moodle/Wiris (ver figura \ref{fig:5_feedback_plataforma}D), la solución paso a paso les genera dudas. No comprenden la factorización que se plantea: $z^3-2= (z-\sqrt[3]{2})(z^2+\sqrt[3]{2}z+\sqrt[3]{2^2})$. Se exhiben algunos fragmentos de diálogos efectuados por los estudiantes en la tabla \ref{tab:4_discusion_factorizacion}. 

\begin{table}[h!]
\caption{Transcripción de la discusión en la clase para intentar comprender el feedback sobre la factorización de $z^3-2$ entregado por la plataforma Moodle/Wiris.}
\label{tab:4_discusion_factorizacion}
\begin{center}
\scalebox{0.7}{\begin{tabular}{p{18cm}l@{}}
\toprule
Paula:   57:50 aquí tengo una duda esto lo elevan todo a tres {[}refiriéndose a la   factorización de $z^3-2$ en $(z-\sqrt{2}   \cdot (z^2-\sqrt[3]{2})z+ \sqrt[3]{4})${]}\\ \midrule
Cristopher: 58:05   es un producto notable ese                                                                    \\ \midrule
Paula: 58:10   es que según yo se puede ocupar elevar al cubo y de ahí nace esto                                                                                              \\ \midrule
Paula: 60:14   quedaría lo mismo que estaba haciendo al principio                                                                                                             \\ \midrule
Paula: 61:21   yo me perdí aquí {[}señala la factorización de una diferencia de cubos que   muestra el feedback de la plataforma{]}                                                                                                                                                                      \\ \midrule
Cristopher: 61:20   la última línea no la entiendo                                                                                                                              \\ \midrule
Juan: 61:51   esa es como una suma por diferencia, pero de cubos o no?                                                                                                                                                                      \\ \midrule
Profesor:   62:06 pensemos en uno más sencillo, que pasa si ustedes intentan en vez de   tener $z$ el cubo es igual a 2 {[}$z^3=2${]} tienen $z$ al cubo igual a 8   {[}$z^3=8${]}. Intenten hacer esas mismas líneas que están ahí {[}refiriéndose al   feedback de Moodle/Wiris{]} y suman el inverso aditivo de 8, les queda $z$ al   cubo menos ocho igual cero {[}$z^3-8=0${]} y si intentan aplicar producto notable   ahí ¿cómo le quedaría lo que le estoy escribiendo en el chat?         \\ \midrule
Cristopher:   66:04 por eso queda en términos de 2 raíz de 3 {[}refiriéndose a $\sqrt[3]{2}${]}   que serían en el caso del 8 sería el 2                                                                                          \\ \midrule
Paula:   66:16 eso le iba a preguntar si yo pusiera $z$ menos 2 elevado al cubo   {[}$(z-2)^3${]} sería lo mismo que decir $z$ al cubo menos 8 {[}$z^3-8${]} y ahí   podríamos hacer la factorización o la suma de cubos                             \\ \midrule
Juan: 67:15   pero es que factorizar es cómo hacer factores y los factores son   los que están en la multiplicación                                                                                                             \\ \midrule
Juan: 67:39   sí, porque convirtió esa expresión a factores, son dos cosas, y la   multiplicación es una operación de 3 factores                                                                                                              \\ \midrule
Cristopher: 68:29   que uno de esos dos factores vale 0                                                           \\ \midrule
Juan:   68:57 queda que $z$ es igual a raíz tercia de 2 {[}$z=\sqrt[3]{2}${]}                                                                                                  \\ \midrule
Paula: 69:55   ¿se puede factorizar la $z$? ¿O no?                                                                                                                               \\ \midrule
Juan: 70:15   a continuación tomamos esa expresión... Ahh no se puede porque   al principio tiene una suma... no no no                                                       \\ \midrule
Paula:   70:43 dejarlo como z y aquí quedaría un 1 y todo esto después se puede pasar   dividendo y se puede racionalizar así se puede utilizar en el $z$   {[}refiriéndose al feedback de Moodle/Wiris{]}                                                                                               \\ \midrule
Paula:   91:23 ya profe, ya intenté factorizar primero las raíces después las $z$ y no   llegué a nada, y me salía error en todo y después pasé esto al otro lado y   tampoco me salía nada. Entonces pase todo al Photomath y salía que había que   resolverlo como una ecuación de segundo grado, da $z_1$ y $z_2$. Y eso que hay que desarrollar la ecuación general y va   a dar el resultado de $z_1$ y $z_2$ que es lo que nos daba en la   solución. \\ \midrule
Paula: 92:28   aquí se desarrolló todo, se aproximó, aquí lo desarrolla todo y lo aproxima,   trabaja con todos los números que se sacaron de la ecuación inicial y ahí las   ecuaciones las aproxima… ve que lo aproxima? Lo trabaja todo con los números   que se sacaron de la ecuación original y al final las ecuaciones las   aproxima.                                                                                                                                              \\ \midrule
Paula:   94:10 ingresé ésta, usted nos pedía desarrollar esta ecuación {[}muestra la   ecuación del feedback (ver figura \ref{fig:5_feedback_plataforma}D){]}                                                                                                             \\ \bottomrule
\end{tabular}}
\end{center}
\end{table}

Consideramos que el siguiente diálogo es rico en la discusión que se genera sobre lo que presenta el feedback, Paula interpreta de forma incorrecta el feedback (ver tabla \ref{tab:4_discusion_factorizacion} minuto 57:50), se da cuenta de su error al dialogar con el profesor. Cristopher reconoce que se usa diferencia de cubos, pero es él mismo, que después indica que no entiende la factorización. Juan luego reconoce algunas nociones, pero con poca seguridad e incluso tienen dudas de qué significa factorizar. El profesor comienza a trabajar en cuestiones más elementales de álgebra, sobre los cuales, los estudiantes, recuerdan haber visto, pero que no dominan, como la factorización de una expresión del tipo $a^3-b^3$. Finalmente, Cristopher reconoce que la expresión $z^3-2$ se puede escribir como diferencia de cubos y factorizar, lo que les permite usar la propiedad “$a \cdot b=0$ implica que $a=0$ o $b=0$” para la solución de la ecuación. En el mismo episodio, Paula está intentando despejar $z$ de una forma que no es posible hacerlo, mostrando debilidades en propiedades algebraicas.  

Se observa que el diálogo sobre factorización toma una gran cantidad de tiempo dentro de la clase. De forma colectiva se va reconstruyendo la factorización.

% Please add the following required packages to your document preamble:
% \usepackage

En términos teóricos, en este episodio observamos un \textit{referencial teórico} del álgebra débil, aunque recuerdan algunos elementos, no logran comprender la factorización propuesta en el feedback Moodle/Wiris. Además hay cuatro dificultades para los estudiantes: la primera es la factorización; [$(z-\sqrt[3]{2}) \cdot (z^2+ \sqrt[3]{2} z+\sqrt[3]{2^2}$]; la segunda es sobre la utilización de la propiedad ``$a \cdot b=0$ implica que $a=0$ o $b=0$'' para descomponer la ecuación en dos ecuaciones más simples; la tercera es sobre la solución de la ecuación de segundo grado que queda al factorizar; y la cuarta fue convertir las soluciones con argumento negativo en la raíz en números complejos. Las primeras tres dificultades están asociadas a conocimientos que se trabajan en la escuela, es decir, dificultades discursivas, en cambio la cuarta dificultad, es un conocimiento nuevo para los estudiantes, pero que no logra ser controlado por el referencial en concstrucción sobre los números complejos.

%corresponden a conocimientos previos  Salvo lo último, el resto son elementos que están de forma superficial en su \textit{referencial teórico}. También se observa que algunos estudiantes no pueden movilizar elementos de control,\footnote{Nos referimos con esto, a que la estudiate al ingresar mal la ecuación, obtiene soluciones reales, pero esto no es coherente con lo que se visualiza en la ecuación en su forma algebraica, dónde tiene dos soluciones complejas y una real} para identificar que la solución entregada por el \textit{artefacto digital} cuando ingresan mal los datos, no es adecuada.

%Sin embargo, cuando Paula lo retoma en el minuto 90 muestra nuevas dificultades, puesto que utiliza Photomath, ingresa uno de los coeficientes de forma errónea (en vez de resolver una ecuación con $\sqrt[3]{4}$ como coeficiente, utilizó $\sqrt{4}$) y no se da cuenta que las soluciones obtenidas no son posibles, puesto que deberían ser complejas. El profesor hace una serie de preguntas sobre lo que están discutiendo, pero los estudiantes parecen con más dudas que al comienzo y es en este momento finaliza la clase. En la clase siguiente, la misma estudiante retoma el problema, indicando que ha modificado el ingreso de coeficientes a Photomath y éste le muestra que la ecuación no tiene una solución en los reales, el profesor guía el proceso de transformación a números complejos.

\section{Conclusiones}

En este artículo se presentaron 2 tareas sobre encontrar la solución de una ecuación en el sistema de los números complejos en un CAA. En éstas se pedían a los estudiantes obtener una de las soluciones de la ecuación $z^n=a$, donde $n$ podía ser 3 o 4. La solución pedida estaba marcada en rojo sobre un cuadrante del plano complejo, e independiente de cómo se obtuvieran las soluciones, los estudiantes debían relacionar el registro algebraico con el gráfico propuesto en la tarea. Se les propuso utilizar distintos CAS a disposición en la web y en celulares para resolver las tareas: Photomath, GeoGebra, Wolfram Alpha, Symbolab y Calcme. A partir de los resultados obtenidos en este estudio, hemos considerado cuatro aspectos relevantes que nos han permitido realizar las conclusiones que se presentan a continuación:
\subsection*{Validez epistemológica de distintos CAS}

%En una primera discusión con el profesor se constató que los estudiantes tenían conocimientos previos sobre los números complejos débiles y aislados. 

Cuando los estudiantes comenzaron a trabajar con los distintos CAS, pudieron constatar que estos sistemas entregan soluciones muy diferentes frente a la misma entrada ``solve $z^n=a$''. Este primer hecho muestra como cada sistema viene con una inteligencia histórica y una validez epistemológica relativa particular \cite{Flores2022}. Cada \textit{artefacto digital} muestra información distinta, ya sea porque las soluciones pueden ser correctas, parcialmente correctas o con algunos errores, o porque utilizan distintos registros de representación semiótica (algebraicos y/o gráficos), creando momentos de oportinudad para mejorar conocimientos de los saberes matemáticos o \textit{instrumentales}.  

%Las dificultades de las respuestas erróneas que puede entregar un \textit{artefacto digital}, si esto no es percibido por el profesor, es difícil que el estudiante pueda darse cuenta. Para sortear esta dificultad es importante comparar la respuesta en distintos artefactos y analizar cuál es la validez de las respuestas que entrega, abriendo la posibilidad a discusiones epistémicamente ricas.

 \subsection*{Potencialidades en el trabajo matemático con artefactos digitales}
Se observaron potencialidades que se pueden explotar para mejorar la comprensión de los estudiantes de los objetos matemáticos involucrados. Uno de ellos fue la articulación de los \textit{artefactos digitales}, uno que entregaba las soluciones de forma algebraica, rectangular y exacta; y otra que permitía una representación gráfica, que dió sentido tanto para comprender la tarea como las soluciones que arrojaban los artecaftos. %Sus compañeros también se beneficiaron de esta discusión.
Al comparar la inteligencia histórica presente en cada artefacto, permite conectar las \textit{génesis semiótica} y \textit{discursiva} para resolver la tarea. Cabe señalar que para superar la imposibilidad de ingresar directamente una raíz no cuadrada, en uno de estos artefactos (GeoGebra), logró movilizar las propiedades de transformar raíces en potencias.

Otras potencialidades observadas fueron la posibilidad que entregan los \textit{artefactos digitales} de guiar el proceso de \textit{justificación}, permitiendo guiar la génesis \textit{discursiva} en el caso de Juan, quien sabiendo los valores aproximados a los que debía llegar, le permitieron construir una \textit{justificación} algebraica de las soluciones. Además, hemos podido constatar que la activación del plano \textit{[Sem-Ins]} dirigió la génesis \textit{discursiva} que no explican los \textit{artefactos digitales}. El estudiante pudo validar sus resultados mediante GeoGebra y Wolfram Alpha, esto nos permite concluir que trabajar con tareas donde se movilicen distintos artefactos puede ser un camino para enfrentar las dificultades que han sido reportadas en distintas investigaciones, como el uso de las técnicas de ``botón'' que invisibilizan procesos matemáticos diversos y fomentan estrategias más superficiales \citep{Lagrange2005a,Jankvist2019}.

Desde el marco teórico utilizado, en esta investigación, podemos señalar que iniciar el trabajo matemático en el plano \textit{[Sem-Ins]}, en la tarea uno, ayudó a los estudiantes a establecer relaciones para identificar la solución pedida en Moodle/Wiris mediante los cuadrantes del plano complejo, y el signo de las soluciones dadas por los \textit{artefactos digitales}. Con esto se observa la activación de la \textit{génesis semiótica} de las representaciones dadas por los artefactos, el registro gráfico y algebraico para identificar la solución pedida. Así %(el caso de Paula en sección 4.2), 
el uso de \textit{artefactos digitales} mejoró la competencia matemática, para pensar en varias representaciones \citep{Caglayan2016}. Cabe señalar, que no todos los estudiantes logran activar estas \textit{génesis} por si mismos, el trabajo \textit{instrumental} en conjunto con el profesor, resultó fundamental para dar sentido a los símbolos entregados por estos \textit{artefactos digitales}. 

%Fue el caso de Sofia, quien al trabajar en paralelo con Symbolab y GeoGebra, y después de superar obstáculos \textit{instrumentales}, los cuales estaban ligados a la \textit{génesis discursiva}  de las raíces (propiedades de potencias), visualiza la solución algebraica en el plano complejo, y logra dar sentido a la solución expresada con radicales dada por Symbolab. 

\subsection*{Dificultades en el trabajo matemático con artefactos digitales}

Se observaron dificultades en la génesis  \textit{instrumentales}, asociadas al ingreso de expresiones en los formatos específicos de cada \textit{artefacto digital}, lo que para algunos estudiantes fue interpretado como un error matemático. También, existen dificultades en la \textit{génesis discursiva}, que se producen por una no compresión entre la información que entrega un \textit{artefacto digital} y el \textit{referencial teórico} de quien recibe e interpreta dicha información. Tal como se observó en el caso de los estudiantes que no lograron comprender los procesos involucrados en la obtención de las soluciones complejas: factorización de una diferencia de cubos, descomposición de la ecuación en ecuaciones de grado menor, solución de la ecuación de grado dos sin soluciones reales y transformación de las soluciones con argumento de la raíz cuadrada negativa en un número complejo. 

En nuestra investigación se pudo evidenciar que el feedback que se les entregó en Moodle/Wiris dio cuenta de dificultades algebraicas de los estudiantes. Una primera dificultad es comprender la propiedad distributiva (factorización); la segunda es sobre la utilización de la propiedad ``$a \cdot b=0$ implica que $a=0$ o $b=0$'' para descomponer la ecuación en dos ecuaciones más simples; la tercera es sobre la solución de la ecuación de segundo grado que queda al factorizar y finalmente es convertir las soluciones con argumento negativo en la raíz en números complejos.

%Todas estas dificultades tienen en común que están asociadas a un \textit{referencial teórico} insuficiente de los estudiantes para superarlas, por lo que el acompañamiento del profesor es indispensable para detectarlas y ayudar a los estudiantes a resolverlas, cuestión, que como se vio en los análisis no es sencilla y requiere que domine no tan solo la matemática involucrada, sino que las diferencias entre los \textit{artefactos digitales}.  

\subsection*{Variables didácticas en el diseño de las tareas y perspectivas}
Finalmente, cobra relevancia la elección de las variables didácticas en el diseño de tareas. En este artículo fue importante que la solución de la tarea dependiera de articular al menos dos registros de representación semiótica. El uso del registro gráfico obliga a trabajar con la \textit{visualización} de las soluciones, promoviendo un tránsito del dominio \citep{MontoyaVivier2014a} del álgebra (ecuaciones) hacia la geometría (polígonos regulares). También, la elección del grado de la ecuación, en este caso, grados 3 y 4, permitió trabajar con aspectos de álgebra como la propiedad distributiva (factorización) y con figuras geométricas tales como el triángulo equilátero y el cuadrado. 

%Con respecto a esto, el plano \textit{[Ins-Dis]} fue activado para verificar si se cumplían las condiciones para determinar que se obtiene un triángulo equilátero o un cuadrado, donde la \textit{génesis discursiva} fue activada al usar la \textit{génesis instrumental}, con el objetivo de comprobar las propiedades de polígonos regulares en el \textit{artefacto digital} GeoGebra. Si se trabaja con una ecuación de grado superior, algunas de estas estrategias dejan de ser posibles o se vuelven más complicadas y se debe trabajar con la fórmula de la raíz n-ésima de un número complejo, por ejemplo. 

Tanto las potencialidades como las limitaciones, muestran la importancia del \textit{referencial teórico} para trabajar con \textit{artefactos digitales}. En el caso de que este sea débil o esté en desarrollo, el trabajo con \textit{artefactos digitales} permite revelar cuáles son los aspectos en el que es necesario mejorar.

Las perspectivas de este trabajo son variadas, por un lado, hacer modificaciones en la tarea, como eliminar el feedback de la solución paso a paso para no inhibir algunas técnicas algebraicas para encontrar las soluciones. Cambiar el grado de la ecuación por uno superior a 4 para que aparezcan otras estrategias que involucren la forma polar o la de Euler, y así profundizar sobre el conocimiento de este sistema numérico. Además, se podría utilizar esta tarea con otro grupo de estudiantes de otro nivel en su formación o en la formación continua de profesores para estudiar la matemática que ellos movilizan y aportar a la construcción de su conocimiento.

En este artículo proponemos como afrontar la creciente oferta de software, que se relaciona con la métofora del Mito de Sísifo \citep{Camus1942}, que obligan a aprender sobre nuevos \textit{artefactos digitales}, confrontando las distintas respuestas que da cada uno y asumiendo que estamos en una época donde los estudiantes tienen a disposición múltiples herramientas, pudiendo usar esto a nuestro favor para fomentar discusiones que enriquezcan el conocimiento matemático de los estudiantes.

%\bibliography{interacttfpsample.bib}
\bibliographystyle{apalike-es.bst}
\bibliography{interacttfpsample}

\end{document}